\documentclass[leqno, 12pt]{article}

\usepackage{amsmath}
\usepackage{amsfonts}
\usepackage{amssymb}
\usepackage{amsthm} 
\usepackage{graphicx}

\usepackage{scalerel,stackengine}
\stackMath
\newcommand\widecheck[1]{%
\savestack{\tmpbox}{\stretchto{%
  \scaleto{%
    \scalerel*[\widthof{\ensuremath{#1}}]{\kern-.6pt\bigwedge\kern-.6pt}%
    {\rule[-\textheight/2]{1ex}{\textheight}}
  }{\textheight}%
}{0.5ex}}%
\stackon[1pt]{#1}{\scalebox{-1}{\tmpbox}}%
}
\usepackage{scalerel,stackengine}
\stackMath
\newcommand\widebreve[1]{%
\savestack{\tmpbox}{\stretchto{%
  \scaleto{%
    \scalerel*[\widthof{\ensuremath{#1}}]{\kern-.1pt\bigcup\kern-.1pt}  
    {\rule[-\textheight/2]{1ex}{\textheight}}
  }{\textheight}%
}{0.5ex}}%
\stackon[1pt]{#1}{\scalebox{1}{\tmpbox}}%
}
\usepackage{comment}

\def\toL1{\buildrel \mathit{L_1}\over\longrightarrow}

\def\towpr1{\buildrel w.pr.1\over\to}
\def\11{\buildrel 1-1\over\longleftrightarrow}

\def\~{\tilde}
 
\def\^{\hat}

\def\implies{\Rightarrow}

\def\sig{\sigma}

\def\Gam{\Gamma}
\def\gam{\gamma}
\def\lam{\lambda}
\def\alp{\alpha}
\def\bet{\beta}
\def\del{\delta}
\def\Del{\Delta}
\def\eps{\epsilon}

\def\Ups{\Upsilon}
\def\var{\vartheta}
\def\ome{\omega}
\def\Ome{\Omega}
\def\kap{\kappa}

\def\prtl{\partial}
\def\nul{\emptyset}
\def\E{{\rm E}}
\def\P{{\rm P}}

\def\R1{{\bf R}^1}
\def\B1{{\cal B}^1}

\def\nid{\noindent}

\def\ts{\textstyle}
\def\t12{\textstyle{1\over2}}

\def\smalltype{\let\rm=\eightrm \let\bf=\eightbf \let\it=\eightit \let\sl=\eightsl
 \baselineskip=9pt \rm}
\font\tenrm=cmr10
\font\tenbf=cmbx10
\font\tenit=cmti10
\font\tensl=cmsl10
\def\medtype{\let\rm=\tenrm \let\bf=\tenbf \let\it=\tenit \let\sl=\tensl
 \baselineskip=12pt \rm}
\font\twelverm=cmr12
\font\twelvebf=cmbx12
\font\twelveit=cmti12 
\font\twelvesl=cmsl12
\def\bigtype{\let\rm=\twelverm \let\bf=\twelvebf \let\it=\twelveit \let\sl=\twelvesl
 \baselineskip=14pt \rm}
 
\begin{document}

\title{Estimating the Ratio of Means in a Zero-inflated Poisson Mixture Model }
\author{Michael D. Perlman\footnote{mdperlma@uw.edu.}\\Department of Statistics\\
University of Washington
}

\maketitle

\begin{abstract}

\nid The problem of estimating the ratio of the means of a two-component Poisson mixture model is considered, when each component is subject to zero-inflation, i.e., excess zero counts.  The resulting {\it zero-inflated Poisson mixture (ZIPM) model} can be treated as a three-component Poisson mixture model with one degenerate component. The EM algorithm is applied to obtain frequentist estimators and their standard errors, the latter determined via an explicit expression for the observed information matrix. Bayes and empirical Bayes estimators also are obtained by means of conjugate priors and their data-based variants. Lastly, the ZIPM distribution and the ZTP (zero-truncated Poisson) distribution are compared.
\end{abstract}

\newpage 

\vskip6pt

\nid{\bf 1. Introduction.} 

Consider an ecological study aimed at determining the relative reproductive rate of a newly discovered invasive subspecies A of ant compared to that of the native subspecies B. The  available data is indirect, consisting only of counts of nests in several standardized sites, rather than direct observations of individuals.  Furthermore, the nests of the two subspecies are indistinguishable, (possibly) differing only in their relative numbers per site. If the expected numbers of nests per site for A and B are denoted by $\mu$ and $\nu$respectively, it is desired to estimate their ratio $\theta\equiv\mu/\nu$, where $0<\theta<\infty$.

Because little is known about the characteristics of  A, no further constraint can be imposed on $\theta$, which renders the problem unidentifiable as stated, i.e. $(\mu,\nu)$ is indistinguishable from $(\nu,\mu)$.
However, it is reasonable to assume that the newly discovered subspecies A is less prevalent than the  established subspecies B, at least initially. This assumption will be incorporated into the mixture model introduced below, rendering it identifiable.

Furthermore, it is typical of such field studies that data is lost due to uncontrollable factors such as rain, resulting in excessive numbers of zero counts.  As is commonly done, we shall adopt the zero-inflated Poisson (ZIP) distribution to represent this feature (cf. Lambert (1992)). 

Let $N_{ij}$ denote the number of ant nests observed on day $i$ at site $j$.  Let ${\cal I}\equiv\{1,\dots,I\}$ and ${\cal J}\equiv\{1,\dots,J\}$ be the corresponding  index sets, and set ${\cal K}={\cal I}\times{\cal J}$, $K=|{\cal K}|=IJ$.
For $(i,j)\in{\cal K}$, consider random variables (rvs)
\begin{align}
Y_j&\sim\mathrm{Bernoulli}(\pi),\label{Ydef}\\
M_{ij}\,|\,Y_j&\sim\mathrm{Poisson}\big\{t_i[0^{Y_j}\mu+(1-0^{Y_j})\nu]\big\};\label{Mdef}\\
N_{ij}&=Z_{ij}M_{ij},\label{Ndef}\\
Z_{ij}&\sim\mathrm{Bernoulli}(\eps);\label{Zdef}
\end{align}
where $0^0=1$, $\{Y_j\}$ and $\{Z_{ij}\}$ are mutually independent, and $\{M_{ij}\}$
and $\{Z_{ij}\}$ are conditionally mutually independent given $\{Y_j\}$. Thus $M_{ij}$ {\it is a $\pi$- mixture} of $\mathrm{Poisson}(t_i\mu)$ and $\mathrm{Poisson}(t_i\nu)$ rvs, where each $t_i>0$ is known, reflecting a daily feature common to all sites, such as temperature,  and $\mu,\nu\in(0,\infty)$ are unknown. Here  $N_{ij}$ is a  {\it zero-inflated Poisson mixture (ZIPM) rv} with zero-inflation parameter $1-\eps\in(0,1)$. 

The main goal of this paper is the problem of estimating the ratio $\theta\equiv\mu/\nu$ based solely on the observed data $\{N_{ij}\}$, with $\{Y_j\}$, $\{M_{ij}\}$, and $\{Z_{ij}\}$ unobserved. As noted above, for identifiability of $(\mu,\nu)$, and therefore of $\theta$, a restriction must be imposed: we assume that $0<\pi\le1/2$, corresponding to the assumption that subspecies A occurs less frequently than subspecies B. Here $\theta,\lam\in(0,\infty)$, where $\lam\equiv\nu$ is viewed as a nuisance parameter. In terms of $(\theta,\lam)$, \eqref{Mdef} can be rewritten as
\begin{align}
M_{ij}\,|\,Y_j\sim\mathrm{Poisson}(t_i\theta^{Y_j}\lam).\label{Mdef2}
\end{align}
Both frequentist and Bayesian analyses will be presented.

Two well-known preliminary problems will serve as guideposts for the main problem. Section 2 reviews the case where ${\{Y_j\}}$ are observed; here inference about $\theta$ is based solely on the Poisson rvs $\{M_{ij}\}$, with $\{Z_{ij}\}$ and $\{N_{ij}\}$ irrelevant. The maximum likelihood estimators (MLEs) $\^\mu$, $\^\nu$, $\^\theta$ and associated confidence intervals are straightforward. For Bayesian analysis (cf. Laurent and Lagrand (2012)), the integrated likelihood function (cf. \eqref{B1})
\begin{align*}
f_\del(\mathbf{m}\,|\,{\bf y};\theta)=\int\nolimits_0^\infty f(\mathbf{m}\,|\,{\bf y};\theta,\lam)\gam_\del(\lam)d\lam
\end{align*}
w.r.to  a gamma prior probability density function (pdf) $\gam_\del(\lam)$ is obtained. A family of conjugate prior pdfs $\phi_{\alp,\bet}(\theta)$ is easily obtained \eqref{D1} from $f_\del(\mathbf{m}\,|\,{\bf y};\theta)$, leading to explicit posterior pdfs, Bayes estimators, and Bayesian credible intervals. Alternatively, the maximum integrated likelihood estimator (MILE), obtained by maximizing $f_\del(\mathbf{m}\,|\,{\bf y};\theta)$ w.r.to $\theta$, is readily determined.

The case where ${\{Y_j\}}$ are unobserved but $\{M_{ij}\}$ are observed is reviewed in Section 3. This can be viewed as a two-component Poisson mixture model for the $\{M_{ij}\}$; again $\{Z_{ij}\}$ and $\{N_{ij}\}$ are irrelevant.  A standard application of the EM algorithm yields the MLEs $\^\pi,\^\mu$, $\^\nu$, and hence $\^\theta$, then their standard errors are approximated via the observed information matrix $I_{\bf m}$, obtained explicitly in \eqref{DDD1}. 

For Bayesian analysis in Section 3, the integrated likelihood function
\begin{align*}
 f_{\var,\del}(\mathbf{y},\mathbf{m}\,|\,\theta)&=\int_0^{1/2}\int_0^\infty  f(\mathbf{y},\mathbf{m}\,|\,\pi,\lam,\theta)\var(\pi)\gam_\del(\lam)d\pi d\lam
\end{align*}
w.r.to $\gam_\del(\lam)$ and any proper prior pdf $\var(\pi)$ for $\pi\in(0,\t12]$ is obtained (cf. \eqref{04A}). From this the integrated likelihood $f_{\var,\del}(\mathbf{m}\,|\,\theta)$ of ${\bf M}$ itself can be found explicitly (cf. \eqref{elemsymm1}). No conjugate prior family is available, but for any prior pdf $\phi(\theta)$ the posterior pdf $f_{\var,\del}(\theta\,|\,\mathbf{m})\propto f_{\var,\del}(\mathbf{m}\,|\,\theta)\phi(\theta)$ can be simulated via MCMC methods, yielding Bayes estimators and credible intervals. 

Alternatively, the conjugate prior $\phi_{\alp,\bet}$ in \eqref{D1} can be replaced by a data-based version that depends on the unobserved $\{Y_j\}$, whose values are then imputed by the EM algorithm, thereby yielding empirical Bayes posterior pdfs, estimators, and credible intervals.

The main problem, where only the ZIPM rvs $\{N_{ij}\}$ are observed, is treated in Section 4. This can be viewed as a three-component Poisson mixture model where one of the components is degenerate at 0.\footnote{A three-component mixture model with two degenerate components, one non-degenerate Poisson component, and i.i.d. observations was considered by Arora and Chaganty  (2021).} Now the EM algorithm yields the MLEs $\^\pi,\eps,\^\mu$, $\^\nu$, and hence $\^\theta$, then their standard errors are approximated via the observed information matrix $I_{\bf n}$, obtained explicitly with some effort in \eqref{00Y}-\eqref{endInfo}, a main contribution of this study. 

For Bayesian analysis in Section 4, the integrated likelihood function
\begin{align}
f_{\var,\eta,\kap,\del}(\mathbf{y},\mathbf{z},\mathbf{n}\,|\,\theta) &=\int_0^{1/2}\int_0^1\int_0^\infty  f(\mathbf{y},\mathbf{z},\mathbf{n}\,|\,\pi,\eps, \theta,\lam)\var(\pi)\xi_{\eta,\kap}(\eps)\gam_\del(\lam)d\pi d\eps d\lam\nonumber 
\end{align}
w.r.to  a gamma prior pdf $\gam_\del(\lam)$, a beta prior $\xi_{\eta,\kap}(\eps)$ for $\eps\in(0,1)$, and any proper prior pdf $\var(\pi)$ for $\pi\in(0,\t12]$  is obtained (cf. \eqref{04B}). From this the integrated likelihoods $f_{\var,\eta,\kap,\del}(\mathbf{z},\mathbf{n}\,|\,\theta)$ and $f_{\var,\eta,\kap,\del}(\mathbf{n}\,|\,\theta)$ can be obtained explicitly, cf. \eqref{elemsymm3} and \eqref{elemsymm4}, although the latter is computationally challenging. Again no conjugate prior family is available, but for any prior pdf $\phi(\theta)$ the posterior pdf $f_{\var,\eta,\kap,\del}(\theta\,|\,\mathbf{n})\propto f_{\var,\eta,\kap,\del}(\mathbf{n}\,|\,\theta)\phi(\theta)$ can be simulated via MCMC methods to obtain Bayes estimators and credible intervals. 

Alternatively, the conjugate prior $\phi_{\alp,\bet}$ in \eqref{D1} can be replaced by a data-based version that depends on the unobserved $\{Y_j\}$ and $\{Z_{ij}\}$, whose values are imputed by the EM algorithm, again yielding empirical Bayes posterior pdfs, estimators, and credible intervals.

The paper concludes with a comparison of the conditional ZIPM distribution and the ZTP (zero-truncated Poisson) distribution in Section 5.
\vskip6pt
\newpage

\nid{\it Notation:} Column vectors and arrays denoted by Roman letters appear in bold type, their components  in plain type; caps denote rvs: 
\begin{align*}
{\bf t}&\equiv(t_1,\dots,t_I)'\in\mathbb{R}^I,\\
{\bf y}&\equiv(y_1,\dots,y_J)'\in\{0,1\}^J,& {\bf Y}&\equiv(Y_1,\dots,Y_J)'\in\{0,1\}^J,\\
{\bf z}&\equiv(z_{ij})\in\{0,1\}^{\cal K},& {\bf Z}&\equiv(Z_{ij})\in\{0,1\}^{\cal K},\\
{\bf m}&=(m_{ij})\in\mathbb{Z}_+^{\cal K},&  {\bf M}&=(M_{ij})\in\mathbb{Z}_+^{\cal K},\\
{\bf n}&=(n_{ij})\in\mathbb{Z}_+^{\cal K},&  {\bf N}&=(N_{ij})\in\mathbb{Z}_+^{\cal K},
\end{align*}
where $\mathbb{Z}_+$ is the set of nonnegative integers. Sums and products will range over the index sets ${\cal I}$ and ${\cal J}$ unless otherwise specified, e.g.,
\begin{align*}
\sum_i&=\sum_{i=1}^I,& \prod_{j}&=\prod_{j=1}^J,& \sum_{i,j}&=\sum_{i=1}^I\sum_{j=1}^J.
\end{align*}
etc. 
Summation over one or neither of the indices $i,j$ involving 
$m_{ij}$, $n_{ij}$, $z_{ij}$, or their random (capitalized) versions will be indicated by simply dropping the indices that are summed over, e.g.,
\begin{align*}
m_i&=\sum_{j}m_{ij},& N_j&=\sum_{i}N_{ij},&\\
m&=\sum_{i,j}m_{ij},& n&=\sum_{i,j}n_{ij}.
\end{align*}

\nid {\bf 2. First preliminary problem: ${\bf Y}$ and ${\bf M}$ observed.} 

\nid Because ${\bf Y}$ is observed, the sets
\begin{align*}
S\equiv S({\bf Y}):&=\{j\in{\cal J}\mid Y_j=1\},\\
T\equiv T({\bf Y}):&=\{j\in{\cal J}\mid Y_j=0\},
\end{align*}
are known, where $S\cup T={\cal J}$. Therefore
the ratio
\begin{align*}
r:\equiv r({\bf Y})&=|S|/J
\end{align*}
also is known, with $1-r=|T|/J$. The conditional probability mass function (pmf) of $\mathbf{M}\equiv\{M_{ij}\}$ given ${\bf Y}$ is 
\begin{align}
f_{\mu,\nu}(\mathbf{m}\,|\,{\bf y})
&=\prod_{i,j}e^{-t_i\theta_j\lam}(t_i\theta_j\lam)^{m_{ij}}/m_{ij}!\label{0}\\
 &=\left(\prod_i\prod_{j\in S}e^{-t_i\mu}\mu^{m_{ij}}\right)\left(\prod_i\prod_{j\in T}e^{-t_i\nu}\nu^{m_{ij}}\right)\cdot\Xi_\mathbf{t}(\mathbf{m})\nonumber\\
 &=e^{-I\bar tJr\mu}\mu^{m_S}\cdot e^{-I\bar tJ(1-r)\nu}\nu^{m_T}\cdot\Xi_\mathbf{t}(\mathbf{m})\nonumber\\
 &=e^{-K\bar tr\mu}\mu^{m_S}\cdot e^{-K\bar t(1-r)\nu}\nu^{m_T}\cdot\Xi_\mathbf{t}(\mathbf{m}),\label{G}
\end{align}
where $\bar t=\frac{1}{I}\sum_it_i$,
\begin{align*}
m_S&=\sum_{j\in S}m_j,\\
m_T&=\sum_{j\in T}m_j,\\
\ts\Xi_\mathbf{t}(\mathbf{m})&\ts=\left(\prod_it_i^{m_{i}}\right)\left(\prod_{i,j}m_{ij}!\right)^{-1}
\end{align*}
Thus, conditional on ${\bf Y}$, $(M_S,M_T)$ is a sufficient statistic for $(\mu,\nu)$, where
\begin{align*}
M_S&\equiv\sum_{j\in S}M_j\sim\mathrm{Poisson}(K\bar tr\mu),\\
M_T&\equiv\sum_{j\in T}M_j\sim\mathrm{Poisson}(K\bar t(1-r)\nu),
\end{align*}
with $M_S$ and $M_T$ independent. 
\vskip4pt

\nid{\bf 2.1. Frequentist analysis.} From \eqref{G}, the MLEs  of $\mu$, $\nu$, and $\theta\equiv\mu/\nu$ are 
\begin{align*}
\^\mu&=\frac{M_S}{K\bar tr },\\
\^\nu&=\frac{M_T}{K\bar t(1-r) },\\
\^\theta&=\frac{(1-r)M_S}{rM_T }.
\end{align*}
Based on $\^\theta$, approximate $(1-\alp)$ confidence intervals for $\theta$ can be developed  in several ways,\footnote{See Li, Tang, and Wong (2014) and references therein.} two of which are presented here. 

First, because $\log\^\theta=\log\^\mu-\log\^\nu$, the familiar normal approximation and propagation of error method shows that for large $K$,
\begin{align*}
\sqrt{K}(\^\mu-\mu)&\approx \ts N(0,\frac{\mu}{\bar tr}),\\
\sqrt{K}(\^\nu-\nu)&\approx \ts N(0,\frac{\nu}{\bar t(1-r)}),\\
\sqrt{K}(\log\^\mu-\log\mu)&\approx \ts N(0,\frac{1}{\bar tr\mu}),\\
\sqrt{K}(\log\^\nu-\log\nu)&\approx \ts N(0,\frac{1}{\bar t(1-r)\nu}),\\
\sqrt{K}(\log\^\theta-\log\theta)&\approx \ts N(0,Q(\bar t;r;\mu,\nu));\\
Q(\bar t;r;\mu,\nu):&=\ts\frac{1}{\bar t}[\frac{1}{r\mu}+\frac{1}{(1-r)\nu}].
\end{align*}
Furthermore, $\^\mu$ and $\^\nu$ are consistent estimators of $\mu$ and $\nu$ respectively, so for sufficiently large $K$, 
\begin{align*}
\sqrt{K}(\log\^\theta-\log\theta)&\approx \ts N(0,Q(\bar t;r;\^\mu,\^\nu)).
\end{align*}
This yields approximate $(1-\alp)$ confidence intervals for $\log\theta$ and $\theta$ given by\footnote{$z_{\alp/2}$ is the upper $(1-\frac{\alp}{2}$-quantile of the standard normal distribution $N(0,1)$.}
\begin{align*}
\log\^\theta\pm\ts \sqrt{\frac{Q(\bar t;r;\^\mu,\^\nu)}{v}}z_{\alp/2},\\
\^\theta\, e^{\pm\sqrt{\frac{Q(\bar t;r;\^\mu,\^\nu)}{v}}z_{\alp/2}},
\end{align*}
respectively, provided that $r$ is not close to 0 or 1.

A second way to obtain an approximate $(1-\alp)$ confidence interval for $\theta$ is to consider the conditional distribution\footnote{Note that this conditional distribution does not depend on $t_1,\dots,t_I$.} of $M_S$ given $M_S+M_T$ (and ${\bf Y}$):
\begin{align*}
[M_S\,|\,M_S+M_T=m]&\sim\mathrm{Binomial}(m,\eta),\\
\eta:=\frac{r\mu}{r\mu+(1-r)\nu}&=\frac{\theta}{\theta+\frac{(1-r)}{r}}\;,
\end{align*}
so $\eta$ is a strictly increasing function of $\theta$. The conditional MLE of $\eta$ is $\^\eta=\frac{M_S}{m}$, so  if $m$ is large\footnote{Note that $\E(m)=\E(M_S+M_T)=K\bar t[r\theta+(1-r)]\lam$.} then 
\begin{align*}
\sqrt{m}(\^\eta-\eta)&\approx N(0,\eta(1-\eta)),\\
\sqrt{m}[2\sin^{-1}(\sqrt{\^\eta})-2\sin^{-1}(\sqrt{\eta})]&\approx N(0,1),
\end{align*}
where $\sin^{-1}(\sqrt{\^\eta})$ is the well-known arcsine variance-stabilizing transformation\footnote{Refinements of this variance-stabilizing transformation are given by Guan (2009).} for the binomial distribution. This yields approximate (conditional) $(1-\alp)$ confidence intervals for $\sin^{-1}(\sqrt{\eta})$ and $\eta$ given by
\begin{align*}
\sin^{-1}(\sqrt{\^\eta})\pm\ts {\frac{1}{2\sqrt{m}}}z_{\alp/2},\\
\left\{\sin\left[\sin^{-1}(\sqrt{\^\eta})\pm\ts {\frac{1}{2\sqrt{m}}}z_{\alp/2}\right]\right\}^2,
\end{align*}
respectively. Because 
\begin{align*}
\theta=\frac{(1-r)}{r}\cdot\frac{1}{\frac{1}{\eta}-1},
\end{align*}
the latter in turn yields an approximate (conditional) $(1-\alp)$ confidence interval for $\theta$ given by
\begin{align*}
\frac{1-r}{r}\cdot\frac{1}{\frac{1}{\left\{\sin\left[\sin^{-1}(\sqrt{\^\eta})\pm\ts {\frac{1}{2\sqrt{m}}}z_{\alp/2}\right]\right\}^2}-1}\;,
\end{align*}
provided that $r$ is not too close to 0 or 1.\footnote{However, Laurent and Lagrand (2012, p.376) state that these conditional confidence intervals may be too conservative.}
\vskip4pt

\nid{\bf 2.2. Bayesian analysis.} Here $\theta$ and $\lam$ are treated as random, so rewrite $f_{\mu,\nu}(\mathbf{m}\,|\,{\bf y})$ in \eqref{0}-\eqref{G} as follows:
\begin{align}
f(\mathbf{m}\,|\,{\bf y};\,\theta,\lam)
&=e^{-K\bar t[r\theta-(1-r)]\lam}\lam^{m_S+m_T}\cdot\theta^{m_S}\cdot\Xi_\mathbf{t}(\mathbf{m}).\label{GG}
\end{align}
The set of gamma pdfs $\gam_\del$ with shape parameter $\del>0$ and scale parameter\footnote{A nonzero scale parameter can be reduced to 1 simply by re-scaling $t_1,\dots,t_J$.} 1 is a conjugate family of prior pdfs for $f(\mathbf{m}\,|\,{\bf y};\,\theta,\lam)$:
\begin{align}
\gam_\del(\lam)=[\Gam(\del)]^{-1}\lam^{\del-1}e^{-\lam},\quad 0<\lam<\infty.\label{gdel}
\end{align}
The integrated pmf of ${\bf M}$ given ${\bf Y}={\bf y}$ is found to be 
\begin{align}
f_\del(\mathbf{m}\,|\,{\bf y};\theta)&=\int\nolimits_0^\infty f(\mathbf{m}\,|\,{\bf y};\theta,\lam)\gam_\del(\lam)d\lam\nonumber\\
&=\frac{\Gam(m+\del)\left(\prod_it_i^{m_i}\right)}{\Gam(\del){\bf m}!}\frac{\theta^{m_S}}{\{K\bar t[r\theta+(1-r)]+1\}^{m_S+m_T+\del} },\label{B1}
\end{align} 
where ${\bf m}!=\prod_{i,j} m_{ij}!$,

It follows from \eqref{B1} that after integrating over $\lam$, $(M_S,M_T)$ is a sufficient statistic for $\theta$ (still requiring that ${\bf Y}$ is known). For fixed $\theta$ and $\lam$,
\begin{align}
[M_S\,|\,Y_j;\theta,\lam]&\sim\mathrm{Poisson}(K\bar tr\theta \lam),\label{MS}\\
[M_T\,|\,Y_j;\theta,\lam]&\sim\mathrm{Poisson}(K\bar t(1-r) \lam),\label{MT}
\end{align}
with $M_S$ and $M_T$ independent, so the joint pmf of $(M_S,M_T)$ given ${\bf Y}={\bf y}$ is
\begin{align*}
f(m_S,m_T\,|\,{\bf y};\theta,\lam)&= e^{-K\bar tr\theta \lam}(K\bar tr\theta \lam)^{m_S}/m_S!\cdot e^{-K\bar t(1-r)\lam}(K\bar t(1-r)\lam)^{m_S}/m_T!\\
&=e^{-K\bar t[r\theta+(1-r)] \lam}\lam^m\cdot (K\bar t)^mr^{m_S}(1-r)^{m_T}/m_S!m_T!
\end{align*}
for $(m_S,m_T)\in\mathbb{Z}_+^2$. Thus the integrated pmf of $(M_S,M_T)$ given ${\bf Y}={\bf y}$ is \begin{align}
f_\del(m_S,m_T\,|\,{\bf y};\theta)&=\int_0^\infty f(m_S,m_T\,|\,{\bf y};\theta,\lam)\gam_\del(\lam)d\lam\nonumber\\
&=\frac{\Gam(m_S+m_T+\del)}{\Gam(\del)m_S!m_T!}\frac{(K\bar t)^{m_S+m_T}r^{m_S}(1-r)^{m_T}\theta^{m_S}}{\{K\bar t[r\theta+(1-r)]+1\}^{m_S+m_T+\del} },\label{MST}
\end{align} 
similar to \eqref{B1}. 

It may interest some to note that \eqref{MST} can be expressed as a generalized bivariate negative binomial pmf:
\begin{align}
f_\del(m_S,m_T\,|\,{\bf y};\theta)
&=\frac{\Gam(m_S+m_T+\del)}{\Gam(\del)m_S!m_T!}a(\theta)^{m_S}b(\theta)^{m_T}(1-a(\theta)-b(\theta))^\del,\label{MST3}\\
a(\theta):&=\frac{K\bar tr\theta}{K\bar tr\theta+K\bar t(1-r)+1},\label{MST4}\\
b(\theta):&=\frac{K\bar t(1-r)}{K\bar tr\theta+K\bar t(1-r)+1}.\label{MST5}
\end{align}

From here there are two paths for inference about $\theta$: 
\vskip4pt

\nid(i) Continue on the Bayesian path and impose a prior distribution\footnote{Under this approach, $\theta$ is independent of $\lam$ a priori, which conforms to the assumption that $\lam$ is a nuisance parameter. Alternatively, some have proposed to treat $\mu$ and $\nu$ as independent a priori and impose separate priors on each. Under the latter approach, however, $\theta=\mu/\nu=\mu/\lam$ so $\theta$ depends on $\lam$, hence $\lam$ is no longer a nuisance parameter.} on $\theta$, from that obtain its posterior distribution; or 

\nid(ii) Carry out frequentist inference about $\theta$ based on \eqref{B1} (or \eqref{MST}).  
\vskip4pt

\nid We shall follow both paths in turn.
\vskip4pt

(i) A conjugate family of prior pdfs for $\theta$ is apparent from \eqref{B1} and \eqref{MST}: 
\begin{align}
\phi_{\alp,\bet;K\bar t,r}(\theta)&=c(\alp,\bet;K\bar t,r)\cdot\frac{\theta^{\alp-1}}{\{K\bar t[r\theta+(1-r)]+1\}^{\alp+\bet}},\quad 0<\theta<\infty;\label{D1}\\
c(\alp,\bet;K\bar t,r))&=\frac{\Gam(\alp+\bet)(K\bar tr)^\alp[K\bar t(1-r)]^\bet}{\Gam(\alp)\Gam(\bet)}.\label{D2}
\end{align}
where $\alp,\bet>0$. These are essentially $F$-densities and require that $r$ is known, which holds  in the present case that ${\bf Y}$ is observed. The prior mean is finite if $\bet>1$ and is given by
\begin{align}
\E_{\alp,\bet;r}[\theta]&=\frac{c(\alp,\bet;K\bar t,r)}{c(\alp+1,\bet-1;K\bar t,r)}\nonumber\\
&=\frac{(1-r)\alp}{r(\bet-1)}.\label{D3}
\end{align}
 From \eqref{B1} and \eqref{D1}, the posterior density of $\theta$ given ${\bf M}={\bf m}$ is 
\begin{align}
f_{\del,\alp,\bet}(\theta\,|\,{\bf y}, {\bf m})&\propto f_\del(\mathbf{m}\,|\,{\bf y};\theta)\phi_{\alp,\bet;K\bar t,r}(\theta)\label{1F}\\
&\propto \frac{\theta^{m_S+\alp-1}}{\{K\bar t[r\theta+(1-r)]+1\}^{m_S+m_T+\alp+\bet+\del} },\label{2F}
\end{align}
hence
\begin{align}
f_{\del,\alp,\bet}(\theta\,|\,{\bf y}, {\bf m})=\phi_{m_s+\alp,m_T+\bet+\del;K\bar t,r}(\theta).\label{4B}
\end{align}
In particular, the Bayes estimator of $\theta$ is given by the posterior mean
\begin{align}
\E_{m_S+\alp,m_T+\bet+\del;r}[\theta\,|\,{\bf y}, {\bf m}]=\frac{(1-r)(m_S+\alp)}{r(m_T+\bet+\del-1)}\label{D4}
\end{align}
if $m_T+\bet+\del>1$. Bayesian posterior confidence intervals for $\theta$ can be obtained from \eqref{4B}.
\vskip4pt

(ii) From \eqref{B1}, after integrating w.r.to $\gam_\del(\lam)$, the conditional log likelihood function of ${\bf M}$ given 
${\bf Y}$ has the form
\begin{align*}
\log f_\del({\bf m}\,|\,{\bf y};\theta)
&=m_S\log\theta-(m_S+m_T+\del)\log\{K\bar t[r\theta+(1-r)]+1\}+h({\bf m}),
\end{align*}
so the conditional maximum integrated likelihood estimator (MILE) of $\theta$ is
\begin{align}
\^\theta_\del({\bf m}\,|\,{\bf y})&=\frac{[K\bar t(1-r)+1]m_S}{K\bar tr(m_T+\del)},\label{E1}
\end{align}
which resembles \eqref{D4}. Theorem 2 of Fahrmeir (1987) applies to show that  if ${\cal I}_{\bf M,\del}(\theta|{\bf y})$ is large\footnote{See condition (D) of Fahrmeir (1987), p.89.}
 then
\begin{align}
\sqrt{K}[\^\theta_\del({\bf m}\,|\,{\bf y})-\theta]\approx N[0,\,K{\cal I}^{-1}_{\bf M,\del}(\theta|{\bf y})],\label{0Q}
\end{align}
where ${\cal I}_{\bf M,\del}(\theta|{\bf y})$ is the expected conditional information number
\begin{align}
{\cal I}_{\bf M,\del}(\theta|{\bf y})&\ts=-\E_\del[\frac{\prtl^2}{\prtl\theta^2}\log f_\del({\bf M}\,|\,{\bf y};\theta)\,|\,{\bf y};\,\theta]\nonumber\\
&=\E_\del\left\{\frac{M_s}{\theta^2}-\frac{(K\bar tr)^2(M_S+M_T+\del)}{\{K\bar t[r\theta+(1-r)]+1\}^2}\bigg|\,{\bf y};\,\theta\right\}\nonumber\\
&=\del K\bar tr\left\{\frac{1}{\theta}-\frac{K\bar tr}{K\bar t[r\theta+(1-r)]+1}\right\}\nonumber\\
&=\frac{\del K\bar tr[K\bar t(1-r)+1]}{\theta\{K\bar t[r\theta+(1-r)]+1\}};\label{1QQ}\\
{\cal I}_{\bf M,\del}^{-1}(\theta|{\bf y})&=\frac{\theta}{\del}\left[\frac{\theta}{K\bar t(1-r)+1}+\frac{1}{K\bar tr}\right].\label{1QQQ}
\end{align}
Here we used the facts that
\begin{align*}
\E_\del(\lam)&=\del;\\
\E_\del[M_S\,|\,{\bf y};\,\theta]&=\E_\del\{\E[M_S\,|\,{\bf y};\,\theta,\,\lam]\,|\,{\bf y};\,\theta\}\\
&=\E_\del\left\{\E\left[\sum\nolimits_i\sum\nolimits_{j\in S}m_{ij}\Big|\,{\bf y};\,\theta,\,\lam\right]\,\Big|\,{\bf y},\,\theta\right\}\\
&=\E_\del\left\{\sum\nolimits_i\sum\nolimits_{j\in S}t_i\theta\lam\,\Big|\,{\bf y};\,\theta\right\}\\
&=K\bar tr\theta\del;\\
\E_\del[M_T\,|\,{\bf y};\,\theta]&=K\bar t(1-r)\del.
\end{align*}

We note from \eqref{1QQ}  that ${\cal I}_{\bf M,\del}(\theta|{\bf y})$ will be large if $K$ is large and $r$ is bounded away from 0 and 1. In this case \eqref{0Q} and \eqref{1QQQ} yield an approximate $(1-\alp)$-confidence interval for $\theta$:
\begin{align}
\^\theta_\del\pm\sqrt{\frac{\^\theta_\del}{\del}\left[\frac{\^\theta_\del}{K\bar t(1-r)+1}+\frac{1}{K\bar tr}\right]}z_{\alp/2}.
\end{align}

\nid {\bf 3. Second preliminary problem: ${\bf Y}$ unobserved, ${\bf M}$ observed.}  

\nid Because $\mathbf{Y}$ is  unobserved, i.e. missing, $S$,  $T$, $r$, $(1-r)$, $M_S$, and $M_T$ are unknown. Here $M_{ij}$ is a $\pi$-mixture of 
$\mathrm{Poisson}(t_i\mu)$ and $\mathrm{Poisson}(t_i\nu)$ rvs, where $\pi$ is the unknown mixing probability, cf. \eqref{Mdef}. 
Thus the (unconditional) pmf of the observed data array $\mathbf{M}\equiv(M_{ij})$ is 
\begin{align}
f_{\pi,\mu,\nu}(\mathbf{m})
&=\prod_{i,j}\left[\pi e^{-t_i\mu}(t_i\mu)^{m_{ij}}+(1-\pi)e^{-t_i\nu}(t_i\nu)^{m_{ij}}\right]/m_{ij}!\label{F}\\\
&=\prod_{i,j}\left[\pi e^{-t_i\mu}\mu^{m_{ij}}+(1-\pi)e^{-t_i\nu}\nu^{m_{ij}}\right]\cdot\Xi_\mathbf{t}(\mathbf{m}).\nonumber
\end{align}
Note that the $K\equiv IJ$ rvs $M_{ij}$ are independent but non-identically distributed (inid) if $t_1,\dots,t_I$ are non-identical. The joint pmf of the complete (unobserved and observed) data $(\mathbf{Y},\mathbf{M})$ is given by
\begin{align}
 f_{\pi,\mu,\nu}(\mathbf{y},\mathbf{m})&=f_\pi({\bf y})f_{\mu,\nu}(\mathbf{m}\,|\,{\bf y})\nonumber\\
 &=\prod_j\pi^{y_j}(1-\pi)^{1-y_j} \cdot\prod_{i,j}\left(e^{-t_i\mu}\mu^{m_{ij}}\right)^{y_j}\left(e^{-t_i\nu}\nu^{m_{ij}}\right)^{1-y_j}\cdot\Xi_\mathbf{t}(\mathbf{m})\label{K}\\
&=\left[\pi^{\bar y}(1-\pi)^{1-\bar y}\right]^J\cdot\prod_j\left[\left(e^{-\bar t\mu}\right)^{I}\mu^{m_j}\right]^{y_j}\left[\left(e^{-\bar t\nu}\right)^{I}\nu^{m_j}\right]^{1-y_j}\cdot\Xi_\mathbf{t}(\mathbf{m}),\label{KKK}\\
&=\left[\pi^{\bar y}(1-\pi)^{1-\bar y}\right]^J\cdot\left[e^{-\bar t\bar y\mu}\mu^{\overline{my}}e^{-\bar t(1-\bar y)\nu}\nu^{\overline{m(1-y)}}\right]^K\cdot\Xi_\mathbf{t}(\mathbf{m}),\label{A}
\end{align}
where $\bar  y=\frac{1}{J}\sum_jy_j$,
\begin{align}
\overline{my}&
=\frac{1}{K}\sum_jm_jy_j,\label{T}\\
\overline{m(1-y)}&=\frac{1}{K}\sum_jm_j(1-y_j)
=\bar m-\overline{my},\label{1TT}\\
\bar m&=\frac{1}{K}\sum_{i,j}m_{ij}=\frac{m}{K}.\label{TT}
\end{align}
Thus $f_{\pi,\mu,\nu}(\mathbf{y},\mathbf{m})$ determines an exponential family with sufficient statistic $(\bar Y, \overline{MY},\overline{M(1-Y)})$, where these are defined similarly to $(\bar y, \overline{my},\overline{m(1-y)})$.
\vskip 4pt

\nid{\it Identifiability:} In the mixture model determined by $f_{\pi,\mu,\nu}$ in \eqref{F}, the parameters $\pi,\mu,\nu$ are not fully identifiable, since $f_{\pi,\mu,\nu}=f_{1-\pi,\nu,\mu}$. Thus, without further specification it is impossible to distinguish between $\theta$ and $1/\theta$. Equivalently, $|\log\theta|$ can be estimated but not $\log\theta$. 

To deal with this, a restriction on the parametrization must be imposed. Often it is assumed that $\mu$ and $\nu$ are ordered, e.g. $\mu<\nu$ which is equivalent to $\theta<1$, but this is inappropriate here. Instead we impose the restriction $\pi\le1/2$, which corresponds to the assumption that $\theta_j=\theta$ occurs less frequently than $\theta_j=1$.
\vskip4pt

\vskip4pt

\nid{\bf 3.1. Frequentist analysis: The EM algorithm.} To obtain the MLEs $\^\pi,\^\mu,\^\nu$ and thus $\^\theta=\^\mu/\^\nu$, it is straightforward to apply the EM algorithm (cf. McLachlan and Krishnan (2008)) as follows:


For $j=1,\dots,J$ define
\begin{align*}
{\bf M}_j&=\{M_{i,j}\,|\,i=1,\dots,I\},\\
{\bf m}_j&=\{m_{i,j}\,|\,i=1,\dots,I\}.
\end{align*}
Because \eqref{A} is an exponential family, for $l=0,1,\dots$, the $(l+1)$-st E-step simply imputes $y_j$ to be
\begin{align}
(\widehat{y_j})_{l+1}&=\E_{\^\pi_l,\^\mu_l,\^\nu_l}[Y_j\,|\,{\bf M}={\bf m}]\nonumber\\
 &=\P_{\^\pi_l,\^\mu_l,\^\nu_l}[Y_j=1\,|\,{\bf M}_j={\bf m}_j]\nonumber\\
 &=\frac{\^\pi_l\prod_{i}e^{-t_i\^\mu_l}(t_i\^\mu_l)^{m_{ij}}}{\^\pi_l\prod_{i}e^{-t_i\^\mu_l}(t_i\^\mu_l)^{m_{ij}}+(1-\^\pi_l)\prod_{i}e^{-t_i\^\nu_l}(t_i\^\nu_l)^{m_{ij}}}\nonumber\\
&=\frac{\^\pi_l}{\^\pi_l+(1-\^\pi_l)e^{-I\bar t(\^\nu_l-\^\mu_l)}(\frac{\^\nu_l}{\^\mu_l})^{m_j}}\label{P}
\end{align}
by Bayes formula. From \eqref{A}, the complete-data MLEs are found to be
\begin{align*}
\~\pi&=\bar  y,\\
\~\mu&=\frac{\overline{my}}{\ \ \bar t\bar y\ \ },\\
\~\nu&=\frac{\overline{m(1-y)}}{\ \ \bar t\,(1-\bar y)\ \ }.
\end{align*}
Thus the $(l+1)$-st M-step yields the updated estimates
\begin{align}
\^\pi_{l+1}&=\frac{1}{J}\sum_j(\widehat{y_j})_{l+1},\label{B}\\
\^\mu_{l+1}&=\frac{\overline{m(\^y)_{l+1}}}{\ \ \bar t\;\overline{(\^y)_{l+1}}\ \ },\label{C}\\
\^\nu_{l+1}&=\frac{\overline{m(1-\^y)_{l+1}}}{\ \ \bar t\,\overline{(1-\^y)_{l+1}}\ \ }
=\frac{\bar m-\overline{m(\^y)_{l+1}}}{\ \ \bar t-\bar t\,\overline{(\^y)_{l+1}}\ \ },\label{D}
\end{align}
where $(\^y)_{l+1}=((\widehat{y_1})_{l+1},\dots,(\widehat{y_J})_{l+1})$. 
If at any stage $\^\pi_{l+1}$ exceeds $1/2$, replace it by $1/2$.

Various improvements to the EM algorithm have been suggested to increase its speed of convergence, etc. See McLachlan and Krishnan (2008) for a thorough survey.

Finally, from \eqref{C} and \eqref{D} we obtain the following updated estimator of $\theta$ (which does not depend on $\bar t$):
\begin{align}
\^\theta_{l+1}&=\frac{\^\mu_{l+1}}{\^\nu_{l+1}}=\frac{\overline{m(\^y)_{l+1}}}{\ \ \overline{(\^y)_{l+1}}\ \ }\frac{1-\overline{(\^y)_{l+1}}}{\ \ \bar m-\overline{m(\^y)_{l+1}}\ \ }=\frac{\frac{1}{\ \overline{(\^y)_{l+1}}\ }-1}{\frac{\bar m}{\ \overline{m(\^y)_{l+1}}\ }-1}\;.\label{E}
\end{align}

\nid{\it Starting value $\^\pi_0$ for the EM algorithm:} Under the restriction $\pi\le1/2$, a simple way to choose $\^\pi_0$ is as follows. Plot a histogram of the entire data set $\{m_{ij}\}$ and attempt to discern two prevalent mixture components, either by eye or by density estimation (cf. Silverman (1986)), then determine their relative weights. Take $\^\pi_0$ to be the lesser of these weights. 
\vskip4pt

\nid{\it Standard error for the MLE $\^\theta$:} For simplicity of notation, set $\ome=(\pi,\mu,\nu)$ and $\^\ome_l=(\^\pi_l,\^\mu_l,\^\nu_l)$. Assume that the EM iterates $\ome_l$ converge to $\^\ome\equiv(\^\pi,\^\mu,\^\nu)$, the actual MLEs based on the observed data ${\bf M}$. Then if $K\equiv IJ$ is large, it follows from Theorem 2 of Hoadley (1971) that
\begin{align}
\sqrt{K}(\^\ome-\ome)\approx N_3[0,\,K{\cal I}^{-1}_{\bf M}(\ome)],\label{Q}
\end{align}
where, with $f_\ome({\bf m})$ given by \eqref{F}, 
\begin{align}
{\cal I}_{\bf M}(\ome)&\equiv-\E_\ome[\nabla_\ome^2\log f_\ome({\bf M})]\label{QQ}
\end{align}
is the total expected information matrix ($3\times3$) for the sample ${\bf M}$.\footnote{For any smooth function $g\equiv g(\ome)$, the gradient $\nabla_\ome g$ is the $3\times1$ vector $(\frac{\prtl g}{\prtl \pi},\frac{\prtl g}{\prtl \mu},\frac{\prtl g}{\prtl \nu})'$ and the Hessian $\nabla_\ome^2g$ is the $3\times3$ matrix of mixed partial derivatives $\frac{\prtl^2 g}{\prtl \pi^2}, \frac{\prtl^2 g}{\prtl \pi\prtl\mu},...,\frac{\prtl^2 g}{\prtl \nu^2}$.} 

However, as noted by Efron and Hinkley (1978) and Louis (1982),  
observed information $I_{\bf m}(\ome)$ usually yields a better normal approximation and often is more readily computed 
than expected information, so we replace \eqref{Q} and \eqref{QQ} by 
\begin{align}
\sqrt{K}(\^\ome-\ome)&\approx N_3[0,\,KI^{-1}_{\bf m}(\ome)],\nonumber\\
I_{\bf m}(\ome)&\equiv-\nabla_\ome^2\log f_\ome({\bf m})\nonumber\\
&=-\E_\ome[\nabla_\ome^2\log f_\ome({\bf m})\,|\,{\bf M}={\bf m}]\nonumber\\
&=-\E_\ome[\nabla_\ome^2\log f_\ome({\bf Y},{\bf m})\,|\,{\bf M}={\bf m}]+\E_\ome[\nabla_\ome^2\log f_\ome({\bf Y}\,|\,{\bf m})\,|\,{\bf M}={\bf m}].\label{Y1}
\end{align}
From \eqref{A},
\begin{align}
 \log f_\ome(\mathbf{Y},\mathbf{m})
 &=J[\bar Y\log\pi+(1-\bar Y)\log(1-\pi)]\nonumber\\
 &\ \ \ \ +K[\overline{mY}\log\mu-\bar t\bar Y\mu+\overline{m(1-Y)}\log\nu-\bar t(1-\bar Y)\nu]+h_\mathbf{t}(\mathbf{m});\nonumber\\ 
\nabla_\ome\log f_\ome({\bf Y},{\bf m})&=\begin{pmatrix}\frac{\ J(\bar Y-\pi)\ }{\pi(1-\pi)}\\ K\left[\frac{\ \overline{mY}\ }{\mu}-\bar t\bar Y\right]\\K\left[\frac{\ \overline{m(1-Y)}\ }{\nu}-\bar t(1-\bar Y)\right]\end{pmatrix};\nonumber\\
-\nabla_\ome^2\log f_\ome({\bf Y},{\bf m})&=\begin{pmatrix}\frac{\ J\left[(1-2\pi)\bar Y+\pi^2\right]\ }{\pi^2(1-\pi)^2}&0&0\\0&K\left[\frac{\ \overline{mY}\ }{\mu^2}\right]&0\\0&0&K\left[\frac{\ \overline{m(1-Y)}\ }{\nu^2}\right]\end{pmatrix};\label{Z}
\end{align}
where $\overline{mY}$ and $\overline{m(1-Y)}$ are defined similarly to \eqref{T} and $h_\mathbf{t}(\mathbf{m})$ does not depend on $\ome$.

Furthermore by \eqref{K}, for fixed ${\bf m}$,
\begin{align}
  f_\ome(\mathbf{y}\,|\,\mathbf{m})&= f_\ome(\mathbf{y},\mathbf{m})/f_{\ome}(\mathbf{m})\nonumber\\
 &\propto\prod_j\left(\pi e^{-\bar tI\mu}\mu^{m_j}\right)^{y_j}\left[(1-\pi )e^{-\bar tI\nu}\nu^{m_j}\right]^{1-y_j},
\end{align}
hence $Y_1,\dots,Y_J$ are conditionally independent given ${\bf M}={\bf m}$ with
\begin{align}
[Y_j\,|\,\mathbf{M}=\mathbf{m}]&\sim\mathrm{Bernoulli}(p_j),\label{W}\\
p_j\equiv p_j(\ome;\bar t;\bar m_j)&=\frac{\pi e^{-\bar tI\mu}\mu^{m_j}}{\pi e^{-\bar tI\mu}\mu^{m_j}+(1-\pi ) e^{-\bar tI\nu}\nu^{m_j}}\;.\nonumber\\
&=\frac{\pi (e^{-\bar t\mu}\mu^{\bar m_j})^{I} }{\pi(e^{-\bar t\mu}\mu^{\bar m_j})^{I} +(1-\pi )(e^{-\bar t\nu}\nu^{\bar m_j})^{I} }\;,\label{WW}
\end{align}
where $\bar m_j=\ts\frac{1}{I}m_j$. From \eqref{W}, 
\begin{align}
\E_\ome[\,\bar Y\,|\,{\bf M}={\bf m}]&=\bar p;\nonumber\\
\E_\ome[\,\overline{mY}\,|\,{\bf M}={\bf m}]&=\overline{mp};\nonumber\\
\E_\ome[\,\overline{m(1-Y)}\,|\,{\bf M}={\bf m}]&=\overline{m(1-p)};\nonumber
\end{align}
where $\overline{mp}$ and $\overline{m(1-p)}$ are defined similarly to \eqref{T}-\eqref{1TT} and 
\begin{align*}
\bar p&=\ts\frac{1}{J}\sum_jp_j,\\
\overline{p(1-p)}&=\ts\frac{1}{J}\sum_jp_j(1-p_j).
\end{align*}
Furthermore, 
\begin{align*}
\E_\ome[\,(1-2\pi)\bar Y+\pi^2\,|\,{\bf M}={\bf m}]&=(1-2\pi)\bar p+\pi^2.
\end{align*}
Thus from \eqref{Z}, the first term in \eqref{Y1} is given by 
\begin{align}
-\E_\ome&[\nabla_\ome^2\log f_\ome({\bf Y},{\bf m})\,|\,{\bf M}={\bf m}]\nonumber\\ \nonumber\\
&=\begin{pmatrix}\frac{J\left[(1-2\pi)\bar p+\pi^2\right]}{\pi^2(1-\pi)^2}&0&0\\0&K\left[\frac{\ \overline{mp}\ }{\mu^2}\right]&0\\0&0&K\left[\frac{\ \overline{m(1-p)}\ }{\nu^2}\right]\end{pmatrix}=:D(\ome;\bar t;{\bf m}).\label{ZZ}
\end{align}

The second term in \eqref{Y1} is obtained as follows: From \eqref{W},
\begin{align*}
 f_\ome(\mathbf{Y}\,|\,\mathbf{m})&=\prod_jp_j^{Y_j}(1-p_j)^{1-Y_j};\\
 \log f_\ome(\mathbf{Y}\,|\,\mathbf{m})&=\sum\limits_jY_j\log p_j+(1-Y_j)\log(1-p_j);\\
 \nabla_\ome\log f_\ome(\mathbf{Y}\,|\,\mathbf{m})&=\sum\limits_j\frac{(Y_j-p_j)}{p_j(1-p_j)}\nabla_\ome p_j;\\
\nabla_\ome^2\log f_\ome(\mathbf{Y}\,|\,\mathbf{m})&\ts=\sum\limits_j\left[\frac{(Y_j-p_j)\nabla_\ome^2p_j-(\nabla_\ome p_j)(\nabla_\ome p_j)'}{p_j(1-p_j)}-\frac{(Y_j-p_j)(\nabla_\ome p_j)[\nabla_\ome(p_j(1-p_j))]'}{p_j^2(1-p_j)^2}\right];\\
\E_\ome[\nabla_\ome^2\log f_\ome(\mathbf{Y}\,|\,\mathbf{m})\,|\,{\bf M}={\bf m}]&=-\sum\limits_j\frac{(\nabla_\ome p_j)(\nabla_\ome p_j)'}{p_j(1-p_j)}\\
&=\sum\limits_j(\nabla_\ome\log p_j)[\nabla_\ome\log(1-p_j)]'.
\end{align*}
From \eqref{WW},
\begin{align*}
\log p_j&=\log \pi-I\bar t \mu+I\bar m_j\log\mu-\log\gam_j,\\
\gam_j\equiv\gam_j(\ome;\bar t;\bar m_j):&=\pi (e^{-\bar t\mu}\mu^{\bar m_j})^{I}+(1-\pi )(e^{-\bar t\nu}\nu^{\bar m_j})^{I},
\end{align*}
from which it can be shown that
\begin{align*}
\nabla_\ome\log p_j&=\frac{(e^{-\bar t\nu}\nu^{\bar m_j})^{I}}{\gam_j}\left(\frac{1}{\pi},(1-\pi)\Big(\frac{\bar m_j}{\mu}-\bar t\Big)I,-(1-\pi)\Big(\frac{\bar m_j}{\nu}-\bar t\Big)I\right)',\\
\nabla_\ome\log(1-p_j)&=-\frac{p_j}{1-p_j}\nabla_\ome\log p_j\\
&=-\frac{\pi}{1-\pi}\Big[e^{\bar t(\nu-\mu)}\Big(\frac{\mu}{\nu}\Big)^{\bar m_j}\Big]^{I}\nabla_\ome\log p_j\\
&=-\frac{(e^{-\bar t\mu}\mu^{\bar m_j})^{I}}{\gam_j}\left(\frac{1}{1-\pi},\pi\Big(\frac{\bar m_j}{\mu}-\bar t\Big)I,-\pi\Big(\frac{\bar m_j}{\mu}-\bar t\Big)I\right)',\\
(\nabla_\ome\log p_j)[\nabla_\ome\log(1-p_j)]'&=-\frac{[e^{-\bar t(\mu+\nu)}(\mu\nu)^{\bar m_j}]^{I}}{\gam_j^2}\del_j\del_j',\\
\del_j\equiv\del_j(\ome;\bar t;m_j):&=\left(\frac{1}{\sqrt{\pi(1-\pi)}},\ I\sqrt{\pi(1-\pi)}\Big(\frac{\bar m_j}{\mu}-\bar t\Big),\ -I\sqrt{\pi(1-\pi)}\Big(\frac{\bar m_j}{\nu}-\bar t\Big)\right)'.
\end{align*}
Therefore
\begin{align}
\E_\ome[\nabla_\ome^2\log f_\ome(\mathbf{Y}\,|\,\mathbf{m})\,|\,{\bf M}={\bf m}]
&=-\sum\limits_j\frac{[e^{-\bar t(\mu+\nu)}(\mu\nu)^{\bar m_j}]^I}{\gam_j^2}\del_j\del_j'.\label{secondterm}
\end{align}
Thus by \eqref{Y1}, \eqref{ZZ},  and \eqref{secondterm}, the observed information matrix is 
\begin{align}
I_{\bf m}(\ome)&=D(\ome;\bar t;{\bf m})-e^{-\bar tI(\mu+\nu)}\Del(\ome;\bar t;{\bf m})\Del(\ome;\bar t;{\bf m})';\label{DDD1}\\
\Del(\ome;\bar t;{\bf m}):&=
\begin{pmatrix}\frac{(\mu\nu)^{I\bar m_1/2}}{\gam_1}\del_1,&\dots\ ,&\frac{(\mu\nu)^{I\bar m_J/2}}{\gam_J}\del_J\end{pmatrix}.\nonumber
\end{align}

 Now estimate $I_{\bf m}(\ome)$ in the normal approximation
 \begin{align}
\sqrt{K}(\^\ome-\ome)&\approx N_3[0,\,KI^{-1}_{\bf m}(\ome)]\nonumber
\end{align}
by replacing $\ome$ in $I_{\bf m}(\ome)$ by its MLE $\^\ome\equiv(\^\pi,\^\mu,\^\nu)$ to obtain
 \begin{align}
\sqrt{K}(\^\ome-\ome)&\approx N_3[0,\,KI^{-1}_{\bf m}(\^\ome)].\label{U}
\end{align}
This requires replacing $\pi,\mu,\nu$ by $\^\pi,\^\mu,\^\nu$ wherever the former three appear in the entries of $I_{\bf m}(\ome)$, including in $p_j$, $\del_j$, and $\gam_j$. For large $K$ the $3\times3$ matrix $I_{\bf m}(\^\ome)$ is positive definite, hence invertible.

Lastly, an approximate confidence interval for $\theta\equiv\mu/\nu\equiv g(\ome)$ is obtained from \eqref{U} via propagation of error: for $\^\theta=\^\mu/\^\nu$,
\begin{align}
\sqrt{K}(\^\theta-\theta)&\approx N[0,\,L(\nabla_\ome g(\ome)|_{\^\ome})'I^{-1}_{\bf m}(\^\ome)\nabla_\ome g(\^\ome)|_{\^\ome}]]\nonumber\\
&= N\left[0,\,K\left(\frac{\prtl g}{\prtl\pi}\Big|_{\^\ome},\frac{\prtl g}{\prtl\mu}\Big|_{\^\ome},\frac{\prtl g}{\prtl\nu}\Big|_{\^\ome}\right)I^{-1}_{\bf m}(\^\ome)\left(\frac{\prtl g}{\prtl\pi}\Big|_{\^\ome},\frac{\prtl g}{\prtl\mu}\Big|_{\^\ome},\frac{\prtl g}{\prtl\nu}\Big|_{\^\ome}\right)'\right]\nonumber\\
&= N\left[0,\,K\left(0,\frac{1}{\^\nu},\frac{-\^\mu}{\^\nu^2}\right)I^{-1}_{\bf m}(\^\ome)\left(0,\frac{1}{\^\nu},\frac{-\^\mu}{\^\nu^2}\right)'\right]\nonumber\\
&=N\left[0,\,K\left(\frac{1}{\^\nu},\frac{-\^\mu}{\^\nu^2}\right)(I_{22}-I_{21}I_{11}^{-1}I_{12})^{-1}\left(\frac{1}{\^\nu},\frac{-\^\mu}{\^\nu^2}\right)'\,\right]\nonumber\\
&\equiv N(0,\^\sig^2).\label{Z4}
\end{align}
where $I_{\bf m}(\^\ome)=\begin{pmatrix}I_{11}&I_{12}\\I_{21}&I_{22}\end{pmatrix}$ with $I_{11}:1\times1$ and $I_{22}:2\times2$. Thus computation of $\^\tau^2$ only requires the inversion of a $2\times2$ matrix. This yields the following approximate $(1-\alp)$ confidence interval for $\theta$:
\begin{align}
\^\theta\pm\frac{\^\sig}{\sqrt{n}}z_{\alp/2}.\label{V}
\end{align}


\nid{\bf 3.2. Bayesian analysis.} Rewrite the joint pmf \eqref{A} of the complete (unobserved and observed) data $(\mathbf{Y},\mathbf{M})$ in terms of $\pi, \theta,\lam$ as follows:
\begin{align}
 f(\mathbf{y},\mathbf{m}\,|\,\pi,\lam,\theta)&=\left[\pi^{\bar y}(1-\pi)^{1-\bar y}\right]^J e^{-\{K\bar t[\bar y\theta+(1-\bar y)]\lam\}}\lam^m\theta^{K\overline{my}}\cdot\Xi_\mathbf{t}(\mathbf{m}),\label{4A}
\end{align} 
since $K[\overline{my}+\overline{m(1-y)}]=K \bar m=m$.
If in addition to the gamma prior density $\gam_\del(\lam)$ for $\lam$ we assume {\it any} proper prior density $\var(\pi)$ for $\pi\in(0,\t12]$, then the integrated joint pmf of $(\mathbf{Y},\mathbf{M})$ is
\begin{align}
 f_{\var,\del}(\mathbf{y},\mathbf{m}\,|\,\theta)&=\int_0^{1/2}\int_0^\infty  f(\mathbf{y},\mathbf{m}\,|\,\pi,\lam,\theta)\var(\pi)\gam_\del(\lam)d\pi d\lam\nonumber\\
 &=
g_\var(J\bar y)\cdot\frac{\Gam(m+\del)\left(\prod_it_i^{m_i}\right)}{\Gam(\del){\bf m}!}\frac{\theta^{K\overline{my}}}{\{K\bar t[\bar y\theta+(1-\bar y)]+1\}^{m+\del}},\label{04A}\\
&\equiv  f_{\var}(\mathbf{y})\cdot  f_{\del}(\mathbf{m}\,|\,\mathbf{y};\theta),\label{04AV}
\end{align}
where $m=\sum\nolimits_{i,j} m_{ij}$ and 
\begin{align*}
g_\var(j)=\int_0^{1/2}\pi^{j}(1-\pi)^{J-j}\var(\pi)d\pi,\quad 0\le j\le J.
\end{align*}

The integrated likelihood $f_{\var,\del}(\mathbf{m}\,|\,\theta)$ of ${\bf M}$ itself can be found explicitly: 
\begin{align}
f_{\var,\del}(\mathbf{m}\,|\,\theta)&=\sum\nolimits_{{\bf y}\in\Ups}f_{\var,\del}(\mathbf{y},\mathbf{m}\,|\,\theta)\nonumber\\
&=\frac{\Gam(m+\del)\left(\prod_it_i^{m_i}\right)}{\Gam(\del){\bf m}!}\sum_{j=0}^J\sum_{\{{\bf y}|J\bar y=j\}}\frac{g_\var(J\bar y)\theta^{K\overline{my}}}{\{K\bar t[\bar y\theta+(1-\bar y)]+1\}^{m+\del}}\nonumber\\
&=\frac{\Gam(m+\del)\left(\prod_it_i^{m_i}\right)}{\Gam(\del){\bf m}!}\sum_{j=0}^J\left\langle\frac{g_\var(j)}{\{K\bar t[(\frac{j}{J})\theta+(1-\frac{j}{J})]+1\}^{m+\del}}\sum_{\{{\bf y}|\sum_{j} y_{j}=j\}}\theta^{\sum_{j}m_{j}y_{j}}\right\rangle\nonumber\\
&=\frac{\Gam(m+\del)\left(\prod_it_i^{m_i}\right)}{\Gam(\del){\bf m}!}\sum_{j=0}^J\frac{g_\var(j)s_j({\bf m};\theta)}{\{K\bar t[(\frac{j}{J})\theta+(1-\frac{j}{J})]+1\}^{m+\del}},\label{elemsymm1}
\end{align}
where $s_j({\bf m};\theta)$ is the $j$-th elementary symmetric function of $\{\theta^{m_j}|j=1,\dots,J\}$:
\begin{align}
s_j({\bf m};\theta)
=\sum_{\substack{\sig\subseteq{\cal J},\\|\sig|=j}}\ \prod_{j\in\sig}\theta^{m_j}=\sum_{\substack{\sig\subseteq{\cal J},\\|\sig|=j}}^{} \theta^{m_\sig},\label{elemsymm2}
\end{align}
$m_\sig=\sum_{j\in\sig}m_j$, and $s_0({\bf m};\theta)=1$.
Thus $(M_1,\dots,M_J)$ is a sufficient statistic for $\theta$ but $f_{\var,\del}(\mathbf{m}\,|\,\theta)$ is not an exponential family, so no conjugate prior is available. However, for any prior density $\phi(\theta)$ the posterior pdf $f_{\var,\del}(\theta\,|\,\mathbf{m})\propto f_{\var,\del}(\mathbf{m}\,|\,\theta)\phi(\theta)$, which is available explicitly via \eqref{elemsymm1}. Thus MCMC methods (cf. Robert and  Casella (2004))
 can be used to obtain the corresponding Bayes estimator and posterior confidence intervals.

Alternatively, it follows from \eqref{04A} and  \eqref{04AV} that
\begin{align}
 f_{\del}(\mathbf{m}\,|\,\mathbf{y};\,\theta)
  &= \frac{\Gam(m+\del)\left(\prod_it_i^{m_i}\right)}{\Gam(\del){\bf m}!}\frac{\theta^{K\overline{my}}}{\{K\bar t[\bar y\theta+(1-\bar y)]+1\}^{m+\del}}.\label{04AX}
\end{align}
Comparing \eqref{04AX} to \eqref{B1} suggests an empirical Bayes approach where $\phi_{\alp,\bet;K\bar t,\bar y}(\theta)$ in \eqref{D1}  is used as a data-based prior density for $\theta$. Here $\bar y$ and $\overline{my}$ are unobserved, but their values can be imputed via the EM algorithm  described above, as follows. 

The EM algorithm will output
\begin{align*}
\^{\bar y}:&=\lim_{l\to\infty}\frac{1}{J}\sum_j(\widehat{y_j})_{l+1}\\
\widehat{\overline{my}}:&=\lim_{l\to\infty}\frac{1}{K}\sum_jm_j(\widehat{y_j})_{l+1},\\
\widehat{\overline{m(1-y)}}:&=\lim_{l\to\infty}\frac{1}{K}\sum_jm_j(1-(\widehat{y_j})_{l+1}),
\end{align*}
(recall \eqref{T}-\eqref{1TT}), where $(\widehat{y_j})_{l+1}$ is given by \eqref{P}-\eqref{D}. Then from \eqref{1F}-\eqref{D4}, noting that $m_S=K\,\overline{my}$ and $m_T=K\,\overline{m(1-y)}$, and replacing $r$ by $\^{\bar y}$, $1-r$ by $1-\^{\bar y}$, $\overline{my}$ by $\widehat{\overline{my}}$, and $\overline{m(1-y)}$ by $\widehat{\overline{m(1-y)}}$, we obtain the empirical Bayes posterior density (compare to \eqref{4B})
\begin{align}
f_{\del,\alp,\bet}(\theta\,|\,\^{\bf y}, {\bf m}):=\phi_{K\widehat{\overline{my}}+\alp,\,K\widehat{\overline{m(1-y)}}+\bet+\del;\,K\bar t,\,\^{\bar y}}(\theta)\label{3B}
\end{align}
and empirical Bayes estimator
\begin{align}
\^\theta_{\del,\alp,\bet}^\mathrm{EB}:=\frac{(1-\^{\bar y})(K\,\widehat{\overline{my}}+\alp)}{\^{\bar y}\left(K\,\widehat{\overline{m(1-y)}}+\bet+\del-1\right)},\label{D5}
\end{align}
provided that $K\,\widehat{\overline{m(1-y)}}+\bet+\del>1$. Empirical Bayes posterior confidence intervals for $\theta$ can be obtained from \eqref{3B}.
\vskip4pt

\nid {\bf Remark 3.1.} Taking $\alp=\bet=0$ leads to the prior density $\phi_{0,0;K\bar t,\bar y}(\theta)=\theta^{-1}$. This is no longer data-based but is an improper prior, hence cannot reflect actual prior knowledge about $\theta$. Nonetheless, proceeding formally from \eqref{3B} and \eqref{D5}, we obtain the posterior density 
\begin{align}
f_{\del,0,0}(\theta\,|\,\^{\bf y}, {\bf m}):=\phi_{K\widehat{\overline{my}},\,K\widehat{\overline{m(1-y)}}+\del;\,K\bar t,\,\^{\bar y}}(\theta),\label{3BB}
\end{align}
which is a proper density if $K\widehat{\overline{my}}>0$,
and from this the empirical Bayes estimator
\begin{align}
\^\theta_{\del,0,0}^\mathrm{EB}:=\frac{(1-\^{\bar y})(K\,\widehat{\overline{my}})}{\^{\bar y}\left(K\,\widehat{\overline{m(1-y)}}+\del-1\right)},\label{D55}
\end{align}
valid if $K\widehat{\overline{m(1-y)}}+\del>1$; this may have desirable frequentist properties.\hfill$\square$
\vskip4pt

\nid {\bf Remark 3.2.} Suppose that we wish to apply the EM algorithm directly to obtain the MILE of $\theta$ based on the integrated joint likelihood $f_{\var,\del}(\mathbf{y},\mathbf{m}\,|\,\theta)$ in \eqref{04A}. Since this is not an exponential family, the E-step is nontrivial, requiring the evaluation of
\begin{align}
G(\theta\,|\,\theta_0;{\bf m}):&=\E_{\var,\del}\left\{\log\left[\frac{f_{\var,\del}(\mathbf{Y}\,|\,\mathbf{m};\theta)}{f_{\var,\del}(\mathbf{Y}\,|\,\mathbf{m};\theta_0)}\right]\,\bigg|\,{\bf m}; \theta_0\right\}\nonumber\\
&=\E_{\var,\del}\left\{\log\left[\frac{f_{\var,\del}(\mathbf{Y},\mathbf{m}\,|\,\mathbf{m};\theta)}{f_{\var,\del}(\mathbf{Y},\mathbf{m}\,|\,\mathbf{m};\theta_0)}\right]\,\bigg|\,{\bf m}; \theta_0\right\}-\log\left[\frac{f_{\var,\del}(\mathbf{m}\,|\,\theta)}{f_{\var,\del}(\mathbf{m}\,|\,\theta_0)}\right]\nonumber\\
&=\E_{\var,\del}\left\{K\,\overline{mY}\log\Big(\frac{\theta}{\theta_0}\Big)-(m+1)\log\left\{\frac{K\bar t[\bar Y\theta+(1-\bar Y)]+1}{K\bar t[\bar Y\theta_0+(1-\bar Y)]+1}\right\}\bigg|\,{\bf m}; \theta_0\right\}\label{0EE}\\
&\ \ \ \ \ -\!\log\left[\sum_{j=0}^J\frac{g_\var(j)s_j({\bf m};\theta)}{\{K\bar t[(\frac{j}{J})\theta\!+\!(1\!-\!\frac{j}{J})]\!+\!1\}^{m+1}}-\sum_{j=0}^J\frac{g_\var(j)s_j({\bf m};\theta_0)}{\{K\bar t[(\frac{j}{J})\theta_0\!+\!(1\!-\!\frac{j}{J})]\!+\!1\}^{m+1}}\right]\nonumber
\end{align}
by \eqref{04A} and \eqref{elemsymm1},
which is to be maximized over $\theta$ in the M-step. However, explicit evaluation of the conditional expectation is problematic. In such cases, approaches such as Monte Carlo simulation have been proposed; cf. McLachlan and Krishnan (2008), Debavelaere and Allassonni\`ere (2021).\hfill$\square$
\vskip6pt

\nid {\bf 4. Main problem: ${\bf Y}$, ${\bf Z}$, ${\bf M}$ unobserved, ${\bf N}$ observed.} 

\nid Here $N_{ij}$ is a zero-inflated Poisson mixture (ZIPM) rv: $N_{ij}$ is an $\eps$-mixture of $M_{ij}$ and $O_{ij}$, where $O_{ij}$ is degenerate at 0, so $O_{ij}\sim\mathrm{Poisson}(\lam=0)$; while $M_{ij}$ is a $\pi$-mixture of 
$\mathrm{Poisson}(t_i\mu)\equiv\mathrm{Poisson}(t_i\theta\lam)$ and $\mathrm{Poisson}(t_i\nu)\equiv\mathrm{Poisson}(t_i\lam)$ rvs. Thus this problem can be viewed as a three-component Poisson mixture model with one degenerate component and non-i.i.d observations. The three weights are
$\pi\eps$, $(1-\pi)\eps$, and $1-\eps$, with the identifiability constraint $0<\pi\le1/2$

For notational simplicity, set $\ome=(\pi,\eps,\mu,\nu)$. Under this three-component mixture model, the unconditional pmf of the observed data $\mathbf{N}\equiv\{N_{ij}\}$ is 
\begin{align*}
f_\ome(\mathbf{n})
&=\prod_{i,j}\left[\pi\eps e^{-t_i\mu}(t_i\mu)^{n_{ij}}+(1-\pi)\eps e^{-t_i\nu}(t_i\nu)^{n_{ij}}+(1-\eps)\,0^{n_{ij}}\right]/{\bf n}!,
\end{align*}
where $0^0=1$ and ${\bf n}!=\prod_{i,j}n_{ij}!$. Again the $K\equiv IJ$ rvs $N_{ij}$ are independent but non-identically distributed (inid) if $t_1,\dots,t_I$ are non-identical. 

The sample space of $({\bf Y}, ({\bf Z},{\bf N}))$ is $\Ups\times\Ome$, where $\Ups=\{0,1\}^{\cal J}$ and
\begin{align*}
\Ome&=\big[\{0,1\}^{\cal K}\times(\mathbb{Z}_+)^{\cal K}\big]\cap\{({\bf z},{\bf n})\,|\,\forall i,j,\, z_{ij}=0\implies n_{ij}=0\}\\
&=\big[\{0,1\}^{\cal K}\times(\mathbb{Z}_+)^{\cal K}\big]\cap\big\{({\bf z},{\bf n})\,\big|\,\forall i,j,\,n_{ij}(1-z_{ij})=0\big\}\\
&=\big[\{0,1\}^{\cal K}\times(\mathbb{Z}_+)^{\cal K}\big]\cap\big\{({\bf z},{\bf n})\,\big|\,\prod_{i,j}0^{n_{ij}(1-z_{ij})}=1\big\},
\end{align*}
with $\mathbb{Z}_+$ the set of nonnegative integers. The joint pmf of the unobserved and observed data $(\mathbf{Y},\mathbf{Z},\mathbf{N})$ on $\Ups\times\Ome$ is
 \begin{align}
 &f_\ome(\mathbf{y},\mathbf{z},\mathbf{n})\nonumber\\
 &=f_\pi({\bf y})f_\eps({\bf z})f_{\mu,\nu}(\mathbf{n}\,|\,{\bf y},\mathbf{z})\nonumber\\
 &=\prod\nolimits_j\pi^{y_j}(1-\pi)^{1-y_j} \prod\nolimits_{i,j}\eps^{z_{ij}}(1-\eps)^{1-z_{ij}}\label{big}\\
 &\ \ \ \cdot\prod\nolimits_{i,j}\left[e^{-t_i\mu}(t_i\mu)^{n_{ij}}\right]^{y_jz_{ij}}\left[e^{-t_i\nu}(t_i\nu)^{n_{ij}}\right]^{(1-y_j)z_{ij}}0^{n_{ij}(1-z_{ij})}\big/{\bf n}!\nonumber\\
&=\left[\pi^{\bar y}(1-\pi)^{1-\bar y}\right]^J\left[\eps^{\bar z}(1-\eps)^{1-\bar z}\right]^K\left[e^{-\overline{tyz}\mu}\mu^{\overline{nyz}}e^{-\overline{t(1-y)z}\,\nu}\nu^{\overline{n(1-y)z}}\,\right]^K\frac{\prod_{i,j}t_i^{n_{ij}z_{ij}}0^{n_{ij}(1-z_{ij})}}{{\bf n}!}\nonumber\\
&=\left[\pi^{\bar y}(1-\pi)^{1-\bar y}\right]^J\left[\eps^{\bar z}(1-\eps)^{1-\bar z}\right]^K\left[e^{-\overline{tyz}\mu}\mu^{\overline{ny}}e^{-\overline{t(1-y)z}\,\nu}\nu^{\overline{n(1-y)}}\,\right]^K\cdot\Xi_\mathbf{t}(\mathbf{z},\mathbf{n}),\label{big2}
\end{align}
where $\mathbf{y}=\{y_j\}$, $\mathbf{z}=\{z_{ij}\}$, $\mathbf{n}=\{n_{ij}\}$,
\begin{align}
\bar  z&=\frac{1}{K}\sum\nolimits_{i,j}z_{ij},\nonumber\\
\overline{tz}&=\frac{1}{K}\sum\nolimits_{i,j}t_{i}z_{ij}=\frac{1}{K}\sum_{i,j}t_{i}z_{ij},\nonumber\\
\overline{tyz}&=\frac{1}{K}\sum\nolimits_{i,j}t_{i}y_jz_{ij}=\frac{1}{K}\sum_{i,j}t_{i}y_jz_{ij},\nonumber
\end{align}
and similarly with $y$ replaced by $1-y$. To obtain \eqref{big2} we have used the facts that for $(\mathbf{y},\mathbf{z},\mathbf{n})\in\Ups\times\Ome$,
\begin{align}
\overline{nyz}&
=\frac{1}{K}\sum\nolimits_{i,j}n_{ij}y_jz_{ij}\nonumber\\
&=\frac{1}{K}\sum\nolimits_{i,j}n_{ij}y_j\nonumber\\
&=\frac{1}{K}\sum\nolimits_jn_jy_j\nonumber\\
&=\overline{ny},\nonumber\\
\Xi_\mathbf{t}(\mathbf{z},\mathbf{n}):&=\frac{\prod_{i,j}t_i^{n_{ij}z_{ij}}0^{n_{ij}(1-z_{ij})}}{{\bf n}!}\nonumber\\
&=\frac{\prod_{i,j}t_i^{n_{ij}}0^{n_{ij}(1-z_{ij})}}{{\bf n}!}\nonumber\\
&=\frac{\prod_it_i^{n_i}\cdot\prod_{i,j}0^{n_{ij}(1-z_{ij})}}{{\bf n}!};\nonumber
\end{align}
and similarly with $y$ replaced by $1-y$. 
Thus $f_{\ome}(\mathbf{y},\mathbf{z},\mathbf{n})$ determines an exponential family with support $\Ups\times\Ome$ and  sufficient statistic
\begin{align*}
(\bar Y, \bar Z,\,\overline{tYZ},\,\overline{t(1-Y)Z},\,\overline{nY},\,\overline{n(1-Y)}).
\end{align*}
\vskip4pt

\nid{\bf 4.1. Frequentist analysis: the EM algorithm.} To obtain the MLEs $\^\eps,\^\pi,\^\mu,\^\nu$ and then $\^\theta=\^\mu/\^\nu$, it is again straightforward - albeit somewhat challenging, including notationally - to apply the EM algorithm, as follows:

For $i=1,\dots,I$ and $j=1,\dots,J$,  define
\begin{align*}
{\bf N}_j&=(N_{i,j}\,|\,i=1,\dots,I),\\
{\bf n}_j&=(n_{i,j}\,|\,i=1,\dots,I);\\
1_{ij}^{\ne}\equiv1_{ij}^{\ne}(n_{ij})&=1-0^{n_{ij}},\\
1_{j}^{\ne}\equiv1_{j}^{\ne}({\bf n}_j)&=\sum\nolimits_i1_{ij}^{\ne},\\
1^{\ne}\equiv1^{\ne}({\bf n})&=\sum\nolimits_j1_{j}^{\ne},\\
t_j^{\ne}\equiv t_j^{\ne}({\bf n}_j)&=\sum\nolimits_i t_i1_{ij}^{\ne};
\end{align*}
\begin{align*}
1_{ij}^=\equiv1_{ij}^=({\bf n}_j)&=0^{n_{ij}},\\
1_{j}^=\equiv1_{j}^=({\bf n}_j)&=\sum\nolimits_i1_{ij}^=,\\
1^=\equiv1^=({\bf n})&=\sum\nolimits_j1_{j}^=,\\
t_j^=\equiv t_j^=({\bf n}_j)&=\sum\nolimits_i t_i1_{ij}^=.
\end{align*}
Here $1_{ij}^{\ne}$ ($1_{ij}^=$) is the indicator function of the event $\{n_{ij}\ne0\}$ ($\{n_{ij}=0\}$), so $1_{j}^{\ne}$ ($1_{j}^=$) is the number of nonzero (zero) $n_{ij}$ with $j$ fixed, etc. Because \eqref{big} is an exponential family, Bayes formula shows that for $l=0,1,\dots$, the $(l+1)$-st E-step imputes $y_j$ as
\begin{align}
(\widehat{y_j})_{l+1}&=\E_{\^\ome_l}[Y_j\,|\,{\bf N}={\bf n}]\nonumber\\
 &=\P_{\^\ome_l}[Y_j=1\,|\,{\bf N}_j={\bf n}_j]\nonumber\\
  &=\P_{\^\ome_l}[Y_j=1]\P_{\^\ome_l}[{\bf N}_j={\bf n}_j\,|\,Y_j=1]/\P_{\^\ome_l}[{\bf N}_j={\bf n}_j]\nonumber\\
 &=\frac{\^\pi_l\prod\limits_{i}[\^\eps_le^{-t_i\^\mu_l}(t_i\^\mu_l)^{n_{ij}}+(1-\^\eps_l)0^{n_{ij}}]}{\^\pi_l\prod\limits_{i}[\^\eps_le^{-t_i\^\mu_l}(t_i\^\mu_l)^{n_{ij}}+(1-\^\eps_l)0^{n_{ij}}]+(1-\^\pi_l)\prod\limits_{i}[\^\eps_le^{-t_i\^\nu_l}(t_i\^\nu_l)^{n_{ij}}+(1-\^\eps_l)0^{n_{ij}}]}\nonumber\\
   &=\frac{A_{j,l}}{A_{j,l}+B_{j,l}};\nonumber
   \end{align}
   where $\^\ome_l=(\^\eps_l,\^\pi_l,\^\mu_l,\^\nu_l)$ (cf. \eqref{B6}-\eqref{zC6}) and 
   \begin{align}
     A_{j,l}:&\ts=\^\pi_l\^\eps_l^{1_{j}^{\ne}}e^{-t_j^{\ne}\^\mu_l}\left(\prod\nolimits_it_i^{n_{ij}}\right)\^\mu_l^{n_{j}}\prod_i[\^\eps_le^{-t_i\^\mu_l}+(1-\^\eps_l)]^{1_{ij}^=},\nonumber\\
 B_{j,l}:&\ts=(1-\^\pi_l)\^\eps_l^{1_{j}^{\ne}}e^{-t_j^{\ne}\^\nu_l}\left(\prod\nolimits_it_i^{n_{ij}}\right)\^\nu_l^{n_{j}}\prod_i[\^\eps_le^{-t_i\^\nu_l}+(1-\^\eps_l)]^{1_{ij}^=}.\nonumber
\end{align}
This simplifies to
\begin{align}
(\widehat{y_j})_{l+1}&=\frac{\^\pi_l}{\^\pi_l+(1-\^\pi_l)e^{-t_j^{\ne}(\^\nu_l-\^\mu_l)}(\frac{\^\nu_l}{\^\mu_l})^{n_j}\prod_i\left[\frac{\^\eps_le^{-t_i\^\nu_l}+(1-\^\eps_l)}{\^\eps_le^{-t_i\^\mu_l}+(1-\^\eps_l)}\right]^{1_{ij}^=}},\label{hatyjell1}
\end{align}
which should be compared to \eqref{P}.

Also at the $(l+1)$-st E-step, $z_{ij}$, $y_jz_{ij}$, and $(1-y_j)z_{ij}$ are imputed as
\begin{align}
(\widehat{z_{ij}})_{l+1}&=\E_{\^\ome_l}[Z_{ij}\,|\,{\bf N}={\bf n}]\nonumber\\
 &=\P_{\^\ome_l}[Z_{ij}=1\,|\,{\bf N}_{ij}={\bf n}_{ij}]\nonumber\\
  &=\P_{\^\ome_l}[Z_{ij}=1]\P_{\^\ome_l}[{\bf N}_{ij}={\bf n}_{ij}\,|\,Z_{ij}=1]/\P_{\^\ome_l}[{\bf N}_{ij}={\bf n}_{ij}]\nonumber\\
 &=\frac{\^\eps_l[\^\pi_le^{-t_i\^\mu_l}(t_i\^\mu_l)^{n_{ij}}+(1-\^\pi_l)e^{-t_i\^\nu_l}(t_i\^\nu_l)^{n_{ij}}]}{\^\eps_l[\^\pi_le^{-t_i\^\mu_l}(t_i\^\mu_l)^{n_{ij}}+(1-\^\pi_l)e^{-t_i\^\nu_l}(t_i\^\nu_l)^{n_{ij}}]+(1-\^\eps_l)0^{n_{ij}}}\nonumber\\
  &=\frac{\^\eps_l[\^\pi_l+(1-\^\pi_l)e^{-t_i(\^\nu_l-\^\mu_l)}(\frac{\^\nu_l}{\^\mu_l})^{n_{ij}}]}{\^\eps_l[\^\pi_l+(1-\^\pi_l)e^{-t_i(\^\nu_l-\^\mu_l)}(\frac{\^\nu_l}{\^\mu_l})^{n_{ij}}]+(1-\^\eps_l)e^{t_i\^\mu_l}(\frac{1}{\^\mu_l})^{n_{ij}}0^{n_{ij}}};\label{hatyjell2}
 \end{align}
  \begin{align}
(\widehat{y_jz_{ij}})_{l+1}&=\E_{\^\ome_l}[Y_jZ_{ij}\,|\,{\bf N}={\bf n}]\nonumber\\
 &=\P_{\^\ome_l}[Y_j=1,Z_{ij}=1\,|\,{\bf N}_j={\bf n}_j]\nonumber\\
  &=\P_{\^\ome_l}[Y_j=1,Z_{ij}=1]\P_{\^\ome_l}[{\bf N}_j={\bf n}_j\,|\,Y_j=1,Z_{ij}=1]/\P_{\^\ome_l}[{\bf N}_j={\bf n}_j]\nonumber\\
 &=\frac{\^\pi_l\^\eps_l\prod\nolimits_{i}e^{-t_i\^\mu_l}(t_i\^\mu_l)^{n_{ij}}}{A_{j,l}+B_{j,l}}\nonumber\\
 &=\frac{\^\pi_l\^\eps_le^{-t_j^=\^\mu_l}\prod_i[\^\eps_le^{-t_i\^\mu_l}+(1-\^\eps_l)]^{-1_{ij}^=}}
  {\^\eps_l^{1_{j}^{\ne}}\left\{\^\pi_l+(1\!-\!\^\pi_l)e^{-t_j^{\ne}(\^\nu_l-\^\mu_l)}(\frac{\^\nu_l}{\^\mu_l})^{n_{j}}\prod_i\left[\frac{\^\eps_le^{-t_i\^\nu_l}+(1-\^\eps_l)}{{\^\eps_le^{-t_i\^\mu_l}+(1-\^\eps_l)}}\right]^{1_{ij}^=} \right\}};\label{hatyjell3}\\
[\widehat{(1-y_j)z_{ij}}]_{l+1}&=(\widehat{z_{ij}})_{l+1} -(\widehat{y_jz_{ij}})_{l+1}.\nonumber
 \end{align}



From \eqref{big2}, the complete-data MLEs are found to be
\begin{align*}
\~\pi&=\bar  y,\\
\~\eps&=\bar z,\\
\~\mu&=\frac{\overline{ny}}{\ \ \overline{tyz}\ \ },\\
\~\nu&=\frac{\overline{n(1-y)}}{\ \ \overline{t(1-y)z}\ \ }.
\end{align*}
Thus the $(l+1)$-st M-step yields the updated estimates
\begin{align}
\^\pi_{l+1}&=\frac{1}{J}\sum_j(\widehat{y_j})_{l+1},\label{B6}\\
\^\eps_{l+1}&=\frac{1}{K}\sum_{i,j}(\widehat{z_{i}})_{l+1},\label{B66}\\
\^\mu_{l+1}&=\frac{\overline{n(\widehat{y})_{l+1}}}{\ \ \overline{t(\widehat{yz})_{l+1}} \ \ }=\frac{\sum_jn_j(\widehat{y_j})_{l+1}}{\ \ \sum_{i,j}t_{i}(\widehat{y_jz_{ij}})_{l+1}\ \ },\label{C6}\\
\^\nu_{l+1}&=\frac{\overline{n[\widehat{(1-y)}]_{l+1}}}{\ \ \overline{t[\widehat{(1-y)z}]_{l+1}} \ \ }=\frac{\sum_{i,j}n_{ij}[1-\widehat{y_j}]_{l+1}}{\ \ \sum_{i,j}t_{i}[\widehat{(1-y_j)z_{ij}}]_{l+1}\ \ }.\label{zC6}
\end{align}
If, at any stage, either $\^\pi_{l+1}$ or $\^\eps_{l+1}$ exceeds $1/2$, replace it by $1/2$.

Finally, the updated estimator $\^\theta_{l+1}\equiv\frac{\^\mu_{l+1}}{\^\nu_{l+1}}$ is obtained from \eqref{C6} and \eqref{zC6}; here, unlike \eqref{E}, it does depend on $\{t_i\}$.
\vskip4pt

\nid{\it Starting value $\^\pi_0$ for the EM algorithm:} Under the restriction $\pi\le1/2$, a simple way to choose $\^\pi_0$ and $\^\eps_0$ is as follows. Plot a histogram of the entire data set $\{n_{ij}\}$ and attempt to discern a spike at 0 and two prevalent mixture components above 0, either by eye or by density estimation, then determine their weights. Take $1-\^\eps_0$ to be the weight of the spike at 0, then take $\^\pi_0$ to be the lesser of the relative weights of the two nonzero components. 
\vskip4pt

\nid{\it Standard error for the MLE $\^\theta$:} Recall that $\ome=(\pi,\eps,\mu,\nu)$ and $\^\ome_l=(\^\pi_l,\^\eps_l,\^\mu_l,\^\nu_l)$ and assume as before that the EM iterates $\ome_l$ converge to $\^\ome\equiv(\^\pi,\^\eps,\^\mu,\^\nu)$, the actual MLEs based on the observed data ${\bf N}$. Again we rely on the results of Hoadley and Efron/Hinkley to provide the normal approximation\footnote{In this section the conditioning events ${\bf Y}={\bf y}$ and ${\bf N}={\bf n}$ are abbreviated as ${\bf y}$ and ${\bf n}$.}
\begin{align}
\sqrt{K}(\^\ome-\ome)&\approx N_4[0,\,KI^{-1}_{\bf n}(\ome)],\label{00YZ}\\
I_{\bf n}(\ome)&\equiv-\nabla_\ome^2\log f_\ome({\bf n})\label{00Y}\\
&=-\E_\ome[\nabla_\ome^2\log f_\ome({\bf n})\,|\,{\bf n}]\nonumber\\
&=-\E_\ome[\nabla_\ome^2\log f_\ome({\bf Y},{\bf Z},{\bf n})\,|\,{\bf n}]+\E_\ome[\nabla_\ome^2\log f_\ome({\bf Y},{\bf Z}\,|\,{\bf n})\,|\,{\bf n}],\label{00Y1}
\end{align}
where $I^{-1}_{\bf n}$ is the $4\times4$ observed information matrix.
If $({\bf Y},{\bf Z},{\bf n})\in\Ups\times\Ome$ then by \eqref{big2}, 
\begin{align}
 \log f_\ome(\mathbf{Y},\mathbf{Z},\mathbf{n})
 &=J[\bar Y\log\pi+(1-\bar Y)\log(1-\pi)]+K[\bar Z\log\eps+(1-\bar Z)\log(1-\eps)]\nonumber\\
 &\ \ \ \ +K[\overline{nY}\log\mu-\overline{tYZ}\mu+\overline{n(1-Y)}\log\nu-\overline{t(1-Y)Z}\nu]+\log\Xi_\mathbf{t}(\mathbf{n},\mathbf{z});\nonumber\\ 
 \nonumber\\
\nabla_\ome\log f_\ome({\bf Y},\mathbf{Z},{\bf n})&=\begin{pmatrix}\frac{\ J(\bar Y-\pi)\ }{\pi(1-\pi)}\\ \frac{\ K(\bar Z-\eps)\ }{\eps(1-\eps)}\\ K\left[\frac{\ \overline{nY}\ }{\mu}-\overline{tYZ}\right]\\K\left[\frac{\ \overline{n(1-Y)}\ }{\nu}-\overline{t(1-Y)Z}\right]\end{pmatrix};\nonumber\\
-\nabla_\ome^2\log f_\ome({\bf Y},\mathbf{Z},{\bf n})&=\begin{pmatrix}\frac{\ J\left[(1-2\pi)\bar Y+\pi^2\right]\ }{\pi^2(1-\pi)^2}&0&0&0\\0&\frac{\ K\left[(1-2\eps)\bar Z+\eps^2\right]\ }{\eps^2(1-\eps)^2} &0&0\\0&0&K\big[\frac{\ \overline{nY}\ }{\mu^2}\big]&0\\0&0&0&K\big[\frac{\ \overline{n(1-Y)}\ }{\nu^2}\big]\end{pmatrix}.\label{ZQ}
\end{align}

Furthermore by \eqref{big}, for $({\bf y},{\bf z},{\bf n})\in\Ups\times\Ome$ with ${\bf n}$ fixed,
\begin{align}
f_\ome&(\mathbf{y},\mathbf{z}\,|\,\mathbf{n})= f_\ome(\mathbf{y},\mathbf{z},\mathbf{n})/f_{\ome}(\mathbf{n})\nonumber\\
& \propto\prod_j\pi^{y_j}(1-\pi)^{1-y_j} \nonumber\\
 &\ \ \cdot\prod_{i,j}\left\{\eps\left[e^{-t_i\mu}(t_i\mu)^{n_{ij}}\right]^{y_j}\left[e^{-t_i\nu}(t_i\nu)^{n_{ij}}\right]^{1-y_j}\right\}^{z_{ij}}[(1-\eps)0^{n_{ij}}]^{1-z_{ij}}.\nonumber
\end{align}
From this, $\{Z_{ij}\}$ are conditionally independent given ${\bf Y}$ and ${\bf N}$, with
\begin{align}
[Z_{ij}\,|\,\mathbf{y},\mathbf{n}]&\sim\mathrm{Bernoulli}(r_{ij}),\label{Zcondindep}\\
r_{ij}\equiv r(t_i;y_j,n_{ij})&\equiv r(\eps,\mu,\nu;t_i;y_j,n_{ij})\nonumber\\
:&=\frac{\eps\left[e^{-t_i\mu}\mu^{n_{ij}}\right]^{y_j}\left[e^{-t_i\nu}\nu^{n_{ij}}\right]^{1-y_j}t_i^{n_{ij}}}{\eps\left[e^{-t_i\mu}\mu^{n_{ij}}\right]^{y_j}\left[e^{-t_i\nu}\nu^{n_{ij}}\right]^{1-y_j}t_i^{n_{ij}}+(1-\eps)0^{n_{ij}}}\nonumber\\
&=1-0^{n_{ij}}+\frac{0^{n_{ij}}\eps\left[e^{-t_i\mu}\mu^{n_{ij}}\right]^{y_j}\left[e^{-t_i\nu}\nu^{n_{ij}}\right]^{1-y_j}t_i^{n_{ij}}}{\eps\left[e^{-t_i\mu}\mu^{n_{ij}}\right]^{y_j}\left[e^{-t_i\nu}\nu^{n_{ij}}\right]^{1-y_j}t_i^{n_{ij}}+(1-\eps)},\label{rijk}
\end{align}
and 
\begin{align*}
&f_{\ome}({\bf y}\,|\,{\bf n})\\
&\propto \sum_{\bf z}f_\ome(\mathbf{y},\mathbf{z}\,|\,\mathbf{n})\\
&\propto
\prod_j\pi^{y_j}(1-\pi)^{1-y_j}\cdot\prod_{i,j}\{\eps\left[e^{-t_i\mu}\mu^{n_{ij}}\right]^{y_j}\left[e^{-t_i\nu}\nu^{n_{ij}}\right]^{1-y_j}t_i^{n_{ij}}+(1-\eps)0^{n_{ij}}\}\\
&=\ \prod_j\pi^{y_j}(1-\pi)^{1-y_j}\\
&\ \ \ \cdot\prod_{\{i,j|n_{ij}\ne0\}}\{\eps\left[e^{-t_i\mu}\mu^{n_{ij}}\right]^{y_j}\left[e^{-t_i\nu}\nu^{n_{ij}}\right]^{1-y_j}t_i^{n_{ij}}\}\cdot\prod_{\{i,j|n_{ij}=0\}}[\eps e^{-t_i\mu y_j}e^{-t_i\nu(1-y_j)}+(1-\eps)]\\
&\propto\ \eps^{1^{\ne}}\prod_j[\pi e^{-t_j^{\ne}\mu}\mu^{n_j}]^{y_j}[(1-\pi)e^{-t_j^{\ne}\nu}\nu^{n_j}]^{1-y_j}\cdot\prod_j\prod_{\{i|n_{ij=0}\}}[\eps e^{-t_i\mu y_j}e^{-t_i\nu(1-y_j)}+(1-\eps)]\\
&\propto\ \prod_j\left\{[\pi e^{-t_j^{\ne}\mu}\mu^{n_j}]^{y_j}[(1-\pi)e^{-t_j^{\ne}\nu}\nu^{n_j}]^{1-y_j}\cdot\prod_i[\eps e^{-t_i\mu y_j}e^{-t_i\nu(1-y_j)}+(1-\eps)]^{1_{ij}^=}\right\}\\
&\propto\prod_jq_j^{y_j}(1-q_j)^{1-y_j},
\end{align*}
where
\begin{align}
q_j&\equiv q_j({\bf n}_j)\equiv q_j(\pi,\eps,\mu,\nu;{\bf t};{\bf n}_j)\nonumber\\
:&=\frac{\pi e^{-t_j^{\ne}\mu}\mu^{n_j}\prod_i[\eps e^{-t_i\mu}+(1-\eps)]^{1_{ij}^=}}{\pi e^{-t_j^{\ne}\mu}\mu^{n_j}\prod_i[\eps e^{-t_i\mu}+(1-\eps)]^{1_{ij}^=}+(1-\pi) e^{-t_j^{\ne}\nu}\nu^{n_j}\prod_i[\eps e^{-t_i\nu}+(1-\eps)]^{1_{ij}^=}}\;.\nonumber
\end{align}
Thus $\{Y_j\}$ are conditionally independent given ${\bf N}$, with
\begin{align}
[Y_j\,|\,\mathbf{n}]&\sim\mathrm{Bernoulli}(q_j).\label{Ycondindep}
\end{align}
Therefore $E_\ome[\,\bar Y\,|\,{\bf n}]=\bar q\equiv \bar q({\bf n})$, while
\begin{align*}
\E_\ome[\,\overline{nY}\,|\,{\bf n}]
&=\frac{1}{K}\sum_{i,j}n_{ij}\E_\ome[Y_j\,|\,{\bf n}]\\
&=\frac{1}{K}\sum_{i,j}n_{ij}q_j\\
&=\overline{nq},
\end{align*}
\begin{align*}
\E_\ome[\,\overline{n(1-Y)}\,|\,{\bf n}]&=\frac{1}{K}\sum_{i,j}n_{ij}(1-q_j)\\
&=\overline{n(1-q)}.
\end{align*}
Furthermore, 
\begin{align}
\E_\ome[\,(1-2\pi)\bar Y+\pi^2\,|\,{\bf n}]&=(1-2\pi)\bar q+\pi^2.\label{neat}
\end{align}

Next,
\begin{align}
\E_\ome[\,\bar Z\,|\,{\bf n}]&=\E_\ome\{\E_\ome[\,\bar Z\,|\,{\bf y},{\bf n}]\,|\,{\bf n}\};\nonumber\\
&=\frac{1}{K}\sum_{i,j}\E_\ome\{r(t_i;Y_j,n_{ij})\,|\,{\bf n}\};\nonumber\\
&=\frac{1}{K}\sum_{i,j}\{q_jr(t_i;1,n_{ij})+(1-q_j)r(t_i;0,n_{ij})\}\label{00A}\\
&\equiv\overline{qr(1)+(1-q)r(0)}.\nonumber
\end{align}
From \eqref{rijk}, note that
\begin{align}
r(t_i;1,n_{ij})&=1-0^{n_{ij}}+\frac{0^{n_{ij}}\eps e^{-t_i\mu}\mu^{n_{ij}}t_i^{n_{ij}}}{\eps e^{-t_i\mu}\mu^{n_{ij}}t_i^{n_{ij}}+(1-\eps)},\label{r1}\\
r(t_i;0,n_{ij})&=1-0^{n_{ij}}+\frac{0^{n_{ij}}\eps e^{-t_i\nu}\nu^{n_{ij}}t_i^{n_{ij}}}{\eps e^{-t_i\nu}\nu^{n_{ij}}t_i^{n_{ij}}+(1-\eps)},\label{r0}
\end{align}
and decompose $\sum_{i,j}$ as $\sum_{\{i,j|n_{ij}\ne0\}}$ + $\sum_{\{i,j|n_{ij}=0\}}$, so \eqref{00A} becomes 
\begin{align}
\E_\ome[\,\bar Z\,|\,{\bf n}]
&=\frac{1}{K}\sum\limits_{i,j}\left[1-0^{n_{ij}}+\frac{q_j0^{n_{ij}}\eps e^{-t_i\mu}\mu^{n_{ij}}t_i^{n_{ij}}}{\eps e^{-t_i\mu}\mu^{n_{ij}}t_i^{n_{ij}}+(1-\eps)}+\frac{(1-q_j)0^{n_{ij}}\eps e^{-t_i\nu}\nu^{n_{ij}}t_i^{n_{ij}}}{\eps e^{-t_i\nu}\nu^{n_{ij}}t_i^{n_{ij}}+(1-\eps)}\right]\nonumber\\
&=\frac{1^{\ne}}{K}+\frac{\eps}{K}\sum\limits_{\{i,j|n_{ij}=0\}}\left[\frac{q_j e^{-t_i\mu}}{\eps e^{-t_i\mu}+(1-\eps)}+\frac{(1-q_j) e^{-t_i\nu}}{\eps e^{-t_i\nu}+(1-\eps)}\right]\nonumber\\
&=\frac{1^{\ne}}{K}+\frac{\eps}{K}\sum\limits_{i,j}
1_{ij}^=\left[\frac{q_j }{\eps+(1-\eps)e^{t_i\mu}}+\frac{1-q_j}{\eps+(1-\eps) e^{t_i\nu}}\right]\nonumber\\
:&=\rho(\eps,\mu,\nu;{\bf t};{\bf n})\equiv\rho.\nonumber
\end{align}
Lastly, 
\begin{align*}
\E_\ome[(1-2\eps)\bar Z+\eps^2\,|\,{\bf n}]&=(1-2\eps)\rho+\eps^2.
\end{align*}

Therefore $-\E_\ome[\nabla_\ome^2\log f_\ome({\bf Y},{\bf Z},{\bf n})\,|\,{\bf n})\,|\,{\bf n}]$ is evaluated explicitly as follows (recall \eqref{00Y1}-\eqref{ZQ}):
\begin{align}
-\E_\ome&[\nabla_\ome^2\log f_\ome({\bf Y},{\bf Z},{\bf n})\,|\,{\bf n})\,|\,{\bf n}]\nonumber\\ &\nonumber\\
&=\begin{pmatrix}\frac{J\left[(1-2\pi)\bar q+\pi^2\right]}{\pi^2(1-\pi)^2}&0&0&0\\0&\frac{\ K\left[(1-2\eps)\rho+\eps^2\right]\ }{\eps^2(1-\eps)^2} &0&0\\0&0&K\left[\frac{\ \overline{nq}\ }{\mu^2}\right]&0\\0&0&0&K\left[\frac{\ \overline{n(1-q)}\ }{\nu^2}\right]\end{pmatrix}.\label{ZZ7}
\end{align}

For the second term in \eqref{00Y1}, it follows from \eqref{Zcondindep} and \eqref{Ycondindep} that 
\begin{align*}
 f_\ome(\mathbf{Y},\mathbf{Z}\,|\,\mathbf{n})&=f_\ome(\mathbf{Y}\,|\,\mathbf{n})f_\ome(\mathbf{Z}\,|\,\mathbf{Y},\mathbf{n})\\
 &=\prod_jq_j^{Y_j}(1-q_j)^{1-Y_j}\prod_{i,}r_{ij}^{Z_{ij}}(1-r_{ij})^{1-Z_{ij}}\\
  &=\prod_jq_j^{Y_j}(1-q_j)^{1-Y_j}\prod_{\{i,j|n_{ij}=0\}}r_{ij}^{Z_{ij}}(1-r_{ij})^{1-Z_{ij}},
  \end{align*}
 since $n_{ij}\ne0\implies Z_{ij}=1$ and $r_{ij}=1$ for $(\mathbf{Y},\mathbf{Z},\mathbf{n})\in\Ups\times\Ome$. Thus
 \begin{align*}
 \log f_\ome(\mathbf{Y},\mathbf{Z}\,|\,\mathbf{n})&=\sum_j[Y_j\log q_j+(1-Y_j)\log(1-q_j)]\\
 &\ \ \ +\sum_{\{i,j|n_{ij}=0\}}[Z_{ij}\log r_{ij}+(1-Z_{ij})\log(1-r_{ij})];\\
 \nabla_\ome\log  f_\ome(\mathbf{Y},\mathbf{Z}\,|\,\mathbf{n})&=\sum\limits_j\frac{(Y_j-q_j)}{q_j(1-q_j)}\nabla_\ome q_j+\sum\limits_{\{i,j|n_{ij}=0\}}\frac{(Z_{ij}-r_{ij})}{r_{ij}(1-r_{ij})}\nabla_\ome r_{ij};\\
\nabla_\ome^2\log  f_\ome(\mathbf{Y},\mathbf{Z}\,|\,\mathbf{n})&\ts=\sum\limits_j\left[\frac{(Y_j-q_j)\nabla_\ome^2q_j-(\nabla_\ome q_j)(\nabla_\ome q_j)'}{q_j(1-q_j)}-\frac{(Y_j-q_j)(\nabla_\ome q_j)[\nabla_\ome(q_j(1-q_j))]'}{q_j^2(1-q_j)^2}\right]\\
+&\ts\sum\limits_{\{i,j|n_{ij}=0\}}\left[\frac{(Z_{ij}-r_{ij})\nabla_\ome^2r_{ij}-(\nabla_\ome r_{ij})(\nabla_\ome r_{ij})'}{r_{ij}(1-r_{ij})}-\frac{(Z_{ij}-r_{ij})(\nabla_\ome r_{ij})[\nabla_\ome(r_{ij}(1-r_{ij}))]'}{r_{ij}^2(1-r_{ij})^2}\right],
\end{align*}
where $r_{ij}\equiv r(t_i;y_j,n_{ij})$. Therefore
\begin{align}
\E_\ome[\nabla_\ome^2\log  f_\ome(\mathbf{Y},\mathbf{Z}\,|\,\mathbf{n})\,|\,{\bf n}]&=-\sum\limits_j\frac{(\nabla_\ome q_j)(\nabla_\ome q_j)'}{q_j(1-q_j)}-\sum\limits_{\{i,j|n_{ij}=0\}}\E_\ome\Big[\frac{(\nabla_\ome r_{ij})(\nabla_\ome r_{ij})'}{r_{ij}(1-r_{ij})}\,\Big|\,{\bf n}\Big]\nonumber\\
&\ts=\sum\limits_j(\nabla_\ome\log q_j)[\nabla_\ome\log(1-q_j)]'\label{twosums}\\
+&\ts\sum\limits_{\{i,j|n_{ij}=0\}}\E_\ome\Big\{(\nabla_\ome\log r(t_i;Y_j,n_{ij]}))[\nabla_\ome\log(1-r(t_i;Y_j,n_{ij}))]'\,\Big|\,{\bf n}\Big\},\nonumber
\end{align}
where we used the facts that for any functions $h({\bf n})$ and $h({\bf y},{\bf n})$,
\begin{align*}
\E_\ome[(Y_j-q_j)h({\bf n})\,|\,{\bf n}]&=h({\bf n})E_\ome[(Y_j-q_j)\,|\,{\bf n}]=0,\\
\E_\ome[(Z_{ij}-r_{ij})h({\bf y},{\bf n})\,|\,{\bf n}]&=\E_\ome\{h({\bf y},{\bf n})\E_\ome[Z_{ij}-r_{ij}\,|\,{\bf y}, {\bf n}]\,|\,{\bf n}\}=0.
\end{align*}

Now note that 
\begin{align*}
\log q_j&=\log \pi-t_j^{\ne}\mu+n_j\log\mu+\sum_i1_{ij}^=\log[\eps(e^{-t_i\mu}-1)+1]-\log\psi_j;\\
\psi_j:&=\pi e^{-t_j^{\ne}\mu}\mu^{n_j}\prod_i[\eps(e^{-t_i\mu}-1)+1]^{1_{ij}^=}+(1-\pi) e^{-t_j^{\ne}\nu}\nu^{n_j}\prod_i[\eps(e^{-t_i\nu}-1)+1]^{1_{ij}^=};\\
\frac{\prtl\psi_j}{\prtl\pi}&= e^{-t_j^{\ne}\mu}\mu^{n_j}\prod\limits_i[\eps(e^{-t_i\mu}-1)+1]^{1_{ij}^=}
-e^{-t_j^{\ne}\nu}\nu^{n_j}\prod\limits_i[\eps(e^{-t_i\nu}-1)+1]^{1_{ij}^=},\\
\frac{\prtl\psi_j}{\prtl\eps}&=\pi e^{-t_j^{\ne}\mu}\mu^{n_j}\sum\limits_i\frac{1_{ij}^=(e^{-t_i\mu}-1)}{\eps(e^{-t_i\mu}-1)+1}\cdot\prod\limits_i[\eps(e^{-t_i\mu}-1)+1]^{1_{ij}^=}
\\
&\ \ \ +(1-\pi)e^{-t_j^{\ne}\nu}\nu^{n_j}\sum\limits_i\frac{1_{ij}^=(e^{-t_i\nu}-1)}{\eps(e^{-t_i\nu}-1)+1}\cdot\prod\limits_i[\eps(e^{-t_i\nu}-1)+1]^{1_{ij}^=},\\
\frac{\prtl\psi_j}{\prtl\mu}&=\pi e^{-t_j^{\ne}\mu}\mu^{n_j}\prod_i[\eps(e^{-t_i\mu}-1)+1]^{1_{ij}^=}\cdot\left[\frac{n_j}{\mu}-t_j^{\ne}-\eps\sum\limits_i\frac{1_{ij}^=t_ie^{-t_i\mu}}{\eps(e^{-t_i\mu}-1)+1}\right],\\
\frac{\prtl\psi_j}{\prtl\nu}&=(1-\pi)e^{-t_j^{\ne}\nu}\nu^{n_j}\prod_i[\eps(e^{-t_i\nu}-1)+1]^{1_{ij}^=}\cdot\left[\frac{n_j}{\nu}-t_j^{\ne}-\eps\sum\limits_i\frac{1_{ij}^=t_ie^{-t_i\nu}}{\eps(e^{-t_i\nu}-1)+1}\right];
\end{align*}
from which it can be shown that
\begin{align*}
\frac{\prtl\log q_j}{\prtl\pi}&=\frac{e^{-t_j^{\ne}\nu}\nu^{n_j}}{\pi\psi_j}\prod\limits_i[\eps(e^{-t_i\nu}-1)+1]^{1_{ij}^=},\\
\frac{\prtl\log q_j}{\prtl\eps}&=\frac{(1-\pi)e^{-t_j^{\ne}\nu}\nu^{n_j}}{\psi_j}\prod_i[\eps(e^{-t_i\nu}-1)+1]^{1_{ij}^=}\sum\limits_i\frac{1_{ij}^=(e^{-t_i\mu}-e^{-t_i\nu})}{[\eps(e^{-t_i\mu}-1)+1][\eps(e^{-t_i\mu}-1)+1]},\\
\frac{\prtl\log q_j}{\prtl\mu}&=\frac{(1-\pi)e^{-t_j^{\ne}\nu}\nu^{n_j}}{\psi_j} \prod_i[\eps(e^{-t_i\nu}-1)+1]^{1_{ij}^=}\cdot\left[\frac{n_j}{\mu}-t_j^{\ne}-\eps\sum\limits_i\frac{1_{ij}^=t_ie^{-t_i\mu}}{\eps(e^{-t_i\mu}-1)+1}\right],\\
\frac{\prtl\log q_j}{\prtl\nu}&=-\frac{(1-\pi)e^{-t_j^{\ne}\nu}\nu^{n_j}}{\psi_j} \prod_i[\eps(e^{-t_i\nu}-1)+1]^{1_{ij}^=}\cdot\left[\frac{n_j}{\nu}-t_j^{\ne}-\eps\sum\limits_i\frac{1_{ij}^=t_ie^{-t_i\nu}}{\eps(e^{-t_i\nu}-1)+1}\right].
\end{align*}
These four partial derivatives determine the $4\times1$ column vector $\nabla_\ome\log q_j$.  Furthermore,
\begin{align*}
\nabla_\ome\log(1-q_j)&=-\frac{q_j}{1-q_j}\nabla_\ome\log q_j\\
&=-\frac{\pi}{1-\pi}\frac{e^{-t_j^{\ne}\mu}\mu^{n_j}\prod_i[\eps(e^{-t_i\mu}-1)+1]^{1_{ij}^=}}{ e^{-t_j^{\ne}\nu}\nu^{n_j}\prod_i[\eps( e^{-t_i\nu}-1)+1]^{1_{ij}^=}} \nabla_\ome\log q_j,
\end{align*}
hence
\begin{align*}
&(\nabla_\ome\log q_j)[\nabla_\ome\log(1-q_j)]'\\
=&-\frac{\pi(1-\pi)e^{-t_j^{\ne}(\mu+\nu)}(\mu\nu)^{n_j}\prod_i\{[\eps(e^{-t_i\mu}-1)+1][\eps(e^{-t_i\nu}-1)+1]\}^{1_{ij}^=}}{\psi_j^2}\phi_j\phi_j',
\end{align*}
where
\begin{align*}
\phi_j\equiv\phi_j(\ome;{\bf t};{\bf n}_j):&=
\begin{pmatrix}
\frac{1}{\pi(1-\pi)}\\ \\\sum\limits_i\frac{1_{ij}^=(e^{-t_i\mu}-e^{-t_i\nu})}{[\eps(e^{-t_i\mu}-1)+1][\eps(e^{-t_i\mu}-1)+1]}\\ \\\left[\frac{n_j}{\mu}-t_j^{\ne}-\eps\sum\limits_i\frac{1_{ij}^=t_ie^{-t_i\mu}}{\eps(e^{-t_i\mu}-1)+1}\right]\\ \\\left[\frac{n_j}{\nu}-t_j^{\ne}-\eps\sum\limits_i\frac{1_{ij}^=t_ie^{-t_i\nu}}{\eps(e^{-t_i\nu}-1)+1}\right]
\end{pmatrix}.
\end{align*}

Next, for $n_{ij}=0$, 
\begin{align*}
r(t_i;1,0)&=\frac{\eps e^{-t_i\mu}}{\eps(e^{-t_i\mu}-1)+1},\\
r(t_i;0,0)&=\frac{\eps e^{-t_i\nu}}{\eps(e^{-t_i\nu}-1)+1};\\
\log r(t_i;1,0)&=\log \eps-t_i\mu-\log[\eps(e^{-t_i\mu}-1)+1],\\
\log(1- r(t_i;1,0))&=\log(1-\eps)-\log[\eps(e^{-t_i\mu}-1)+1];\\
\log r(t_i;0,0)&=\log \eps-t_i\nu-\log[\eps(e^{-t_i\nu}-1)+1],\\
\log(1- r(t_i;0,0))&=\log(1-\eps)-\log[\eps(e^{-t_i\nu}-1)+1];
\end{align*}
so with $\ome=(\pi,\eps,\mu,\nu)$,
\begin{align*}
\nabla_\ome\log r(t_i;1,0)&=\left(0,\;\frac{1}{\eps[\eps(e^{-t_i\mu}-1)+1]},\;\frac{-(1-\eps)t_i}{\eps(e^{-t_i\mu}-1)+1},\;0\right)',\\
\nabla_\ome\log(1- r(t_i;1,0))&=\left(0,\;\frac{-e^{-t_i\mu}}{(1-\eps)[\eps(e^{-t_i\mu}-1)+1]},\;\frac{\eps t_ie^{-t_i\mu}}{\eps(e^{-t_i\mu}-1)+1},\;0\right)',\\
\nabla_\ome\log r(t_i;0,0)&=\left(0,\;\frac{1}{\eps[\eps(e^{-t_i\nu}-1)+1]},\;0,\;\frac{-(1-\eps)t_i}{\eps(e^{-t_i\nu}-1)+1}\right)',\\
\nabla_\ome\log(1- r(t_i;0,0))&=\left(0,\;\frac{-e^{-t_i\nu}}{(1-\eps)[\eps(e^{-t_i\nu}-1)+1]},\;0,\;\frac{\eps t_ie^{-t_i\nu}}{\eps(e^{-t_i\nu}-1)+1}\right)'.
\end{align*}
Thus 
\begin{align*}
&\E_\ome\Big\{(\nabla_\ome\log r(t_i;Y_j,0))[\nabla_\ome\log(1-r(t_i;Y_j,0))]'\,\Big|\,{\bf n}\Big\}\\
=&\ \ \ (\nabla_\ome\log r(t_i;1,0))[\nabla_\ome\log(1-r(t_i;1,0))]'q_j\\
&\hskip0pt+(\nabla_\ome\log r(t_i;0,0))[\nabla_\ome\log(1-r(t_i;0,0))]'(1-q_j)\}\\
=&-\frac{\eps(1-\eps)e^{-t_i\mu}}{[\eps(e^{-t_i\mu}-1)+1]^2}\,\chi_i(1)\chi_i(1)'q_j\\
 &-\frac{\eps(1-\eps)e^{-t_i\nu}}{[\eps(e^{-t_i\nu}-1)+1]^2}\,\chi_i(0)\chi_i(0)'(1-q_j),
\end{align*}
where
\begin{align*}
\chi_i(1)\equiv\chi(\eps;t_i;1):=\left(0,\ \frac{1}{\eps(1-\eps)},\ t_i,\ 0\right)',\\
\chi_i(0)\equiv\chi(\eps;t_i;0):=\left(0,\ \frac{1}{\eps(1-\eps)},\ 0,\ t_i\right)'.
\end{align*}
Therefore from \eqref{twosums}, the second term in \eqref{00Y1} is given by 
\begin{align}
&\E_\ome[\nabla_\ome^2\log  f_\ome(\mathbf{Y},\mathbf{Z}\,|\,\mathbf{n})\,|\,{\bf n}]\label{endInfo}\\
=&-\pi(1-\pi)\sum_j\frac{e^{-t_j^{\ne}(\mu+\nu)}(\mu\nu)^{n_j}\prod_i\{[\eps(e^{-t_i\mu}-1)+1][\eps(e^{-t_i\nu}-1)+1]\}^{1_{ij}^=}}{\psi_j^2}\phi_j\phi_j'\nonumber\\
&-\eps(1-\eps)\sum\limits_{i,j}1_{ij}^=\bigg\{\frac{e^{-t_i\mu}}{[\eps(e^{-t_i\mu}-1)+1]^2}\,\chi_i(1)\chi_i(1)'q_j+\frac{e^{-t_i\nu}}{[\eps(e^{-t_i\nu}-1)+1]^2}\,\chi_i(0)\chi_i(0)'(1-q_j)\bigg\}.\nonumber
\end{align}
{\it Together with \eqref{ZZ7}, this explicitly determines the observed information matrix $I_{\bf n}(\ome)$ in \eqref{00YZ}-\eqref{00Y1}.}


 Now estimate $I_{\bf n}(\ome)$ in the normal approximation
 \begin{align}
\sqrt{K}(\^\ome-\ome)&\approx N_4[0,\,KI^{-1}_{\bf n}(\ome)]\nonumber
\end{align}
by replacing $\ome$ in $I_{\bf n}(\ome)$ by its MLE $\^\ome\equiv(\^\pi,\^\eps,\^\mu,\^\nu)$, obtained via the EM algorithm,  to obtain
 \begin{align}
\sqrt{K}(\^\ome-\ome)&\approx N_4[0,\,KI^{-1}_{\bf n}(\^\ome)],\label{UA}
\end{align}
where $K=IJ$. This requires replacing $\pi,\eps,\mu,\nu$ by $\^\pi,\^\eps,\^\mu,\^\nu$ wherever the former appear in the entries of $I_{\bf n}(\^\ome)$, including in $q_j$, $\rho$, $\psi_j$, $\phi_j$, and $\chi_i$. For large $K$ the $4\times4$ matrix  $I_{\bf n}(\^\ome)$ is positive definite, hence invertible.

Lastly, an approximate confidence interval for $\theta\equiv\mu/\nu\equiv g(\ome)$ is obtained from \eqref{UA} by propagation of error. For $\^\theta=\^\mu/\^\nu$,
\begin{align}
\sqrt{K}(\^\theta-\theta)&\approx N[0,\,K(\nabla_\ome g(\ome)|_{\^\ome}])'I^{-1}_{\bf n}(\^\ome)\nabla_\ome g(\^\ome)|_{\^\ome}]\nonumber\\
&=N\left[0,\,K\left(\frac{\prtl g}{\prtl\pi}\Big|_{\^\ome},\frac{\prtl g}{\prtl\eps}\Big|_{\^\ome},\frac{\prtl g}{\prtl\mu}\Big|_{\^\ome},\frac{\prtl g}{\prtl\nu}\Big|_{\^\ome}\right)I^{-1}_{\bf n}(\^\ome)\left(\frac{\prtl g}{\prtl\pi}\Big|_{\^\ome},\frac{\prtl g}{\prtl\eps}\Big|_{\^\ome},\frac{\prtl g}{\prtl\mu}\Big|_{\^\ome},\frac{\prtl g}{\prtl\nu}\Big|_{\^\ome}\right)'\,\right]\nonumber\\
&=N\left[0,\,K\left(0,0,\frac{1}{\^\nu},\frac{-\^\mu}{\^\nu^2}\right)I^{-1}_{\bf n}(\^\ome)\left(0,0,\frac{1}{\^\nu},\frac{-\^\mu}{\^\nu^2}\right)'\,\right]\nonumber\\
&=N\left[0,\,K\left(\frac{1}{\^\nu},\frac{-\^\mu}{\^\nu^2}\right)(I_{22}-I_{21}I_{11}^{-1}I_{12})^{-1}\left(\frac{1}{\^\nu},\frac{-\^\mu}{\^\nu^2}\right)'\,\right]\nonumber\\
&\equiv N(0,\^\tau^2),\label{Z41}
\end{align}
where $I_{\bf n}(\^\ome)=\begin{pmatrix}I_{11}&I_{12}\\I_{21}&I_{22}\end{pmatrix}$ is the partitioning of 
$I_{\bf n}(\^\ome)$ into $2\times2$ blocks. Thus computation of $\^\tau^2$ only requires the inversion of two $2\times2$ matrices. This yields the following approximate $(1-\alp)$ confidence interval for $\theta$:
\begin{align}
\^\theta\pm\frac{\^\tau}{\sqrt{K}}z_{\alp/2}.\label{VD}
\end{align}

\nid {\bf 4.2. Bayesian analysis.} Rewrite the joint pmf \eqref{big2} of the complete (unobserved and observed) data $(\mathbf{Y},\mathbf{Z},\mathbf{N})$ in terms of the  parameters $\pi,\eps,\theta,\lam$ as follows:  
\begin{align}
 f(\mathbf{y},\mathbf{z},\mathbf{n}\,|\,\pi,\eps,\theta,\lam)
 &=\left[\pi^{\bar y}(1-\pi)^{1-\bar y}\right]^J\left[\eps^{\bar z}(1-\eps)^{1-\bar z}\right]^K\nonumber\\
 &\ \ \ \cdot\left[e^{-\overline{tyz}\theta\lam}(\theta\lam)^{\overline{ny}}e^{-\overline{t(1-y)z}\,\lam}\lam^{\overline{n(1-y)}}\,\right]^K\cdot\Xi_\mathbf{t}(\mathbf{z},\mathbf{n})\nonumber\\
 &=\left[\pi^{\bar y}(1-\pi)^{1-\bar y}\right]^J\left[\eps^{\bar z}(1-\eps)^{1-\bar z}\right]^K\label{big3}\\
 &\ \ \ \cdot e^{-K[\overline{tyz}\theta+\overline{t(1-y)z}]\,\lam}\lam^{n}\cdot\theta^{K\overline{ny}}\cdot\Xi_\mathbf{t}(\mathbf{z},\mathbf{n}),\nonumber
\end{align}
where $n=\sum_{i,j}n_{ij}$. If we assume the gamma prior density $\gam_\del(\lam)$ for $\lam$, {\it any} proper prior density $\var(\pi)$ for $\pi\in(0,\t12)$, and the beta$(\eta,\kap)$ prior density
\begin{align*}
\xi_{\eta,\kap}(\eps):=\frac{\Gam(\eta+\kap)}{\Gam(\eta)\Gam(\kap)}\eps^{\eta-1}(1-\eps)^{\kap-1}1_{(0,1)}(\eps)
\end{align*}
for $\eps$, where $\eta,\kap>0$, then from \eqref{big3} the integrated joint pmf of $(\mathbf{Y}, \mathbf{Z},\mathbf{N})$ on $\Ups\times\Ome$ is 
\begin{align}
f_{\var,\eta,\kap,\del}(\mathbf{y},\mathbf{z},\mathbf{n}\,|\,\theta) &=\int_0^{1/2}\int_0^1\int_0^\infty  f(\mathbf{y},\mathbf{z},\mathbf{n}\,|\,\pi,\eps, \theta,\lam)\var(\pi)\xi_{\eta,\kap}(\eps)\gam_\del(\lam)d\pi d\eps d\lam\nonumber\\
 &=g_\var(J\bar y)h_{\eta,\kap}(K\bar z)\frac{\Gam(n+\del)\left(\prod_it_i^{n_i}\right)}{\Gam(\del){\bf n}!}\frac{\theta^{v\overline{ny}}}{\{K[\,\overline{tyz}\theta+\overline{t(1-y)z}\,]+1\}^{n+\del}} 1_\Ups({\bf y})1_\Ome(\mathbf{z},\mathbf{n}),\label{04B}
 \end{align}
 (compare to \eqref{04A}), where for $0\le j\le J$ and $0\le \ell\le K$,
 \begin{align}
g_\var(j)&=\int_0^{1/2}\pi^{j}(1-\pi)^{J-j}\var(\pi)d\pi,\nonumber\\
h_{\eta,\kap}(\ell)&=\frac{\Gam(\eta+\kap)\Gam(\eta+\ell)\Gam[\kap+K-\ell]}{\Gam(\eta)\Gam(\kap)\Gam(\eta+\kap+K)}.\nonumber
\end{align}

For $\sig\subseteq{\cal J}$ and $\tau\subseteq{\cal K}$ let $1_\sig$ and $1_\tau$ denote their indicator functions. From \eqref{04B}, the integrated joint pmf $f_{\var,\eta,\kap,\del}(\mathbf{z},\mathbf{n}\,|\,\theta)$ of $({\bf Z},\mathbf{N})$  can be expressed explicitly as follows: 
\begin{align}
&\ f_{\var,\eta,\kap,\del}(\mathbf{z},\mathbf{n}\,|\,\theta)\nonumber\\
=&\sum\nolimits_{{\bf y}\in\Ups}f_{\var,\eta,\kap,\del}(\mathbf{y},\mathbf{z},\mathbf{n}\,|\,\theta)\nonumber\\
=&\frac{\Gam(n+\del)\left(\prod_it_i^{n_i}\right)}{\Gam(\del){\bf n}!}\sum_{j=0}^J\sum_{\{{\bf y}|J\bar y=j\}}\frac{g_\var(J\bar y)h_{\eta,\kap}(K\bar z)\theta^{K\overline{ny}}}{\{K[\,\overline{tyz}\theta+\overline{t(1-y)z}\,]+1\}^{n+\del}}1_\Ome(\mathbf{z},\mathbf{n})\nonumber\\
=&\frac{\Gam(n+\del)\left(\prod_it_i^{n_i}\right)}{\Gam(\del){\bf n}!}\sum_{j=0}^Jg_\var(j)\left\langle\sum_{\sig\subseteq{\cal J},\,|\sig|=j}\frac{\theta^{n_\sig}}{\{(tz)_\sig(\theta-1)+ (tz)+1\}^{n+\del}}\right\rangle h_{\eta,\kap}(K\bar z)1_\Ome(\mathbf{z},\mathbf{n}),\label{elemsymm3}
\end{align}
where
\begin{align*}
n_\sig&=\sum_jn_j1_\sig(j),\\
(tz)_\sig&=\sum_{i,j}t_i1_\sig(j)z_{ij}\\
(tz)&=\sum_{i,j}t_iz_{ij}=K\,\overline{tz},
\end{align*}
and $n_\nul=(tz)_\nul=0$. 

Furthermore, the integrated likelihood $f_{\var,\eta,\kap,\del}(\mathbf{n}\,|\,\theta)$ of $\mathbf{N}$ itself can be obtained explicitly from \eqref{04B} as follows. Setting $\Ome_{\bf n}=\{{\bf z}|({\bf z},{\bf n})\in\Ome\}$,
\begin{align}
&\ f_{\var,\eta,\kap,\del}(\mathbf{n}\,|\,\theta)\nonumber\\
=&\sum\nolimits_{{\bf y}\in\Ups}\sum\nolimits_{{\bf z}\in\Ome_{\bf n}}f_{\var,\eta,\kap,\del}(\mathbf{y},\mathbf{z},\mathbf{n}\,|\,\theta)\nonumber\\
=&\frac{\Gam(n+\del)\left(\prod_it_i^{n_i}\right)}{\Gam(\del){\bf n}!}\sum_{j=0}^J\sum_{\ell=0}^K\sum_{\{{\bf y}|J\bar y=j\}}\sum_{\{{\bf z}|v\bar z=\ell\}}\frac{g_\var(J\bar y)h_{\eta,\kap}(K\bar z)\theta^{K\overline{ny}}1_{\Ome_{\bf n}}({\bf z})}{\{K[\,\overline{tyz}\theta+\overline{t(1-y)z}\,]+1\}^{n+\del}}\nonumber\\
&=\frac{\Gam(n+\del)\left(\prod_it_i^{n_i}\right)}{\Gam(\del){\bf n}!}\sum\limits_{j=0}^J\sum\limits_{\ell=0}^Kg_\var(j)h_{\eta,\kap}(\ell)\left\langle\sum\limits_{\sig\subseteq{\cal J},\,|\sig|=j}\theta^{n_\sig}\Del_{\ell,\sig}({\bf n}\,|\,\theta)\right\rangle,\label{elemsymm4}
\end{align}
where
\begin{align}
n_\sig&=\sum_jn_j1_\sig(j),\label{nsig}\\
\Del_{\ell,\sig}({\bf n}\,|\,\theta)&=\sum\limits_{\tau\subseteq{\cal K},\,|\tau|=\ell}\frac{1_{\Ome_{\bf n}}(1_\tau)}{\{t_{\sig,\tau}(\theta-1)+ t_\tau+1\}^{n+\del}},\label{Del}\\
t_{\sig,\tau}&=\sum_{i,j}t_i1_\sig(j)1_\tau(i,j),\nonumber\\
t_\tau&=\sum_{i,j}t_i1_\tau(i,j),\nonumber
\end{align}
and $t_{\nul,\tau}=t_{\sig,\nul}=t_\nul=0$. 

Because $f_{\var,\eta,\kap,\del}(\mathbf{n}\,|\,\theta)$ is not an exponential family, no conjugate prior is available. However, for any prior density $\phi(\theta)$ the posterior pdf
\begin{align*}
f_{\var,\eta,\kap,\del}(\theta\,|\,\mathbf{n})\propto f_{\var,\eta,\kap,\del}(\mathbf{n}\,|\,\theta)\phi(\theta),
\end{align*}
 which can be obtained  explicitly\footnote{In principle. There are a total of $2^J$ subsets $\sig\subseteq{\cal J}$ and $2^K$ subsets $\tau\subseteq{\cal K}$ that appear in the summations in \eqref{elemsymm4}, where $K=IJ$, so exact calculation of $f_{\var,\eta,\kap,\del}(\mathbf{n}\,|\,\theta)$ is infeasible if $K$ is large. Instead, Monte Carlo simulation over $(\sig,\tau)$ can be used to approximate $f_{\var,\eta,\kap,\del}(\mathbf{n}\,|\,\theta)$} via \eqref{elemsymm4}-\eqref{Del}.  Thus MCMC methods (Robert and Casella (2004)) can be used to simulate the posterior distribution of $\theta$ and thereby obtain the corresponding Bayes estimator and posterior confidence intervals.

Alternatively, we can adopt an empirical Bayes approach as in Section 3.2.  For the data-based prior pdf $\phi_{\alp,\bet;K\overline{tz},\bar r}(\theta)$ (cf. \eqref{D1}), where $\bar r=\frac{\overline{tyz}}{\overline{tz}}$, it follows from \eqref{04B} and \eqref{D1} that the integrated posterior pdf of $\theta$, given the complete data $({\bf y}, ({\bf z},{\bf n}))\in\Ups\times\Ome$,  satisfies
\begin{align}
f_{\var,\eta,\kap,\del}(\theta\,|\,\mathbf{y},\mathbf{z},\mathbf{n})
&\propto\ f_{\var,\eta,\kap,\del}(\mathbf{y},\mathbf{z},\mathbf{n}\,|\,\theta)\phi_{\alp,\bet;K\overline{tz},\bar r}(\theta)\nonumber\\
&\propto\ \frac{\theta^{K\overline{ny}+\alp-1}}{\{K[\,\overline{tyz}\theta+\overline{t(1-y)z}\,]+1\}^{n+\alp+\bet+\del}} 
\nonumber\\
&\propto \phi_{K\overline{ny}+\alp,K\overline{n(1-y)}+\bet+\del;K\overline{tz},\bar r}(\theta),\label{04AY}
\end{align}
since $n=K\overline{ny}+K\overline{n(1-y)}$; note that \eqref{04AY} does not depend on $\var,\eta,\kap$.
Here $K\overline{ny}$, $K\overline{n(1-y)}$, $K\overline{tz}$, $K\overline{tyz}$, and thus $\bar r$, are unobserved, but we can impute their values via the above-discussed EM algorithm as follows:

The EM algorithm will output
\begin{align}
K\widehat{\overline{ny}}&=\lim_{l\to\infty}\sum_jn_j(\widehat{y_j})_{l+1},\label{aa}\\
K\widehat{\overline{n(1-y)}}&=n-K\widehat{\overline{ny}},\label{bb}\\
K\widehat{\overline{tz}}&=\lim_{l\to\infty}\sum_{i,j}t_i(\widehat{z_{ij}})_{l+1},\label{cc}\\
K\widehat{\overline{tyz}}&=\lim_{l\to\infty}\sum_{i,j}t_i(\widehat{y_jz_{ij}})_{l+1},\label{dd}\\
\^{\bar r}&=\frac{\widehat{\overline{tyz}}}{\widehat{\overline{tz}}},\label{ee}
\end{align}
where $(\widehat{y_j})_{l+1}$, $(\widehat{z_{ij}})_{l+1}$, and $(\widehat{y_jz_{ij}})_{l+1}$ appear in \eqref{hatyjell1}-\eqref{hatyjell3}. Now refer to \eqref{1F}-\eqref{D4} and  replace  $r$ by $\^{\bar r}$, $m_S$ by $K\widehat{\overline{ny}}$, and $m_T$ by $K\widehat{\overline{n(1-y)}}$, thus we obtain the empirical Bayes integrated posterior density 
\begin{align}
f_{\del,\alp,\bet}(\theta\,|\,\^{\bf y}, \^{\bf z}, {\bf n}):=\phi_{K\widehat{\overline{ny}}+\alp,\,K\widehat{\overline{n(1-y)}}+\bet+\del;\,K\widehat{\overline{tz}},\,\^{\bar r}}(\theta)\label{3BW}
\end{align}
and empirical Bayes estimator\footnote{Note that the imputed values of $\{z_{ij}\}$ occur in $\^\theta_{\del,\alp,\bet}^\mathrm{EBZIP}$ through $K\widehat{\overline{ny}}$ and $K\widehat{\overline{n(1-y)}}$ as well as through $\^{\bar r}$; see  \eqref{aa}-\eqref{bb}, \eqref{hatyjell1}, and \eqref{B6}-\eqref{zC6}.}
\begin{align}
\^\theta_{\del,\alp,\bet}^\mathrm{EBZIP}:=\frac{(1-\^{\bar r})(K\,\widehat{\overline{ny}}+\alp)}{\^{\bar r}\left(K\,\widehat{\overline{n(1-y)}}+\bet+\del-1\right)},\label{D5W}
\end{align}
provided that $K\,\widehat{\overline{n(1-y)}}+\bet+\del>1$. Empirical Bayes integrated posterior confidence intervals for $\theta$ can be obtained from \eqref{3BW}.
\vskip4pt

\nid {\bf Remark 4.1.} Taking $\alp=\bet=0$ yields the prior density $\phi_{0,0;\,K\overline{tz},\^{\bar r}}(\theta)=\theta^{-1}$. This is no longer data-based but is improper, hence cannot reflect actual prior knowledge about $\theta$. However, proceeding formally from \eqref{3BW} and \eqref{D5W}, we obtain the posterior density 
\begin{align}
f_{\del,0,0}(\theta\,|\,\^{\bf y}, \^{\bf z},{\bf n}):=\phi_{K\widehat{\overline{ny}},\,K\widehat{\overline{n(1-y)}}+\del;\,K\widehat{\overline{ny}},\,\^{\bar r}}(\theta),\label{3BBG}
\end{align}
which is a proper density if $K\widehat{\overline{ny}}>0$, 
and from this the estimator
\begin{align}
\^\theta_{\del,0,0}^{\mathrm{EBZIP}}:=\frac{(1-\^{\bar r})(K\,\widehat{\overline{ny}})}{\^{\bar r}\left(K\,\widehat{\overline{n(1-y)}}+\del-1\right)},\label{D556}
\end{align}
valid if $K\widehat{\overline{n(1-y)}}+\del>1$, and which may have desirable frequentist properties.\hfill$\square$
\vskip4pt

\nid{\bf Remark 4.2.} Direct determination of the MILE of $\theta$ based on $f_{\var,\eta,\kap,\del}(\mathbf{n}\,|\,\theta)$ in \eqref {elemsymm4} again appears problematic. As in Remark 3.2, one might attempt to obtain this MILE by applying the EM algorithm to $f_{\var,\eta,\kap,\del}(\mathbf{y},\mathbf{z},\mathbf{n}\,|\,\theta)$ in \eqref{04B} or to $f_{\var,\eta,\kap,\del}(\mathbf{z},\mathbf{n}\,|\,\theta)$ in \eqref{elemsymm3}, but again the E-steps are challenging.\hfill$\square$
\vskip4pt

\nid{\bf Remark 4.3.} Note that the term $1_{\Ome_{\bf n}}(1_\tau)$ in $\Del_{\ell,\sig}({\bf n}\,|\,\theta)$ (cf. \eqref{Del}) depends on ${\bf n}$ only through
\begin{align*}
1-0^{\bf n}:&=(1-0^{n_{ij}}\,|\,(i,j)\in{\cal K})\\
&=\{1_{ij}^{\ne}\,|\,(i,j)\in{\cal K}),
\end{align*}
i.e., the indicator function over ${\cal K}\equiv{\cal I}\times{\cal J}$ of the set of nonzero $n_{ij}$'s. Thus we obtain the following interesting fact from \eqref{elemsymm4}-\eqref{Del} and the Factorization Criterion: $(N_1,\dots,N_J;1-0^{\bf N})$ is a sufficient statistic\footnote{This holds for any choice of the prior pdf $\var(\pi)$.}
for $\theta$ based on the integrated likelihood $f_{\var,\eta,\kap,\del}(\mathbf{n}\,|\,\theta)$. If we recall that in the non-ZIP model of Part I, $(M_1,\dots,M_J)$ is a sufficient statistic for $\theta$ based on the integrated likelihood $f_{\var,\del}(\mathbf{m}\,|\,\theta)$ for $\theta$ given in \eqref{elemsymm1}, then this shows that in the Bayesian framework, after integrating over the parameters $\pi,\eps,\lam$, the statistic $1-0^{\bf N}$ is the only additional information needed for inference about $\theta$ when zero-inflation is present. This raises the interesting  question of determining the joint distribution of $(N_1,\dots,N_J;1-0^{\bf N})$ based on the integrated likelihood $f_{\var,\eta,\kap,\del}(\mathbf{n}\,|\,\theta)$.\hfill$\square$
\vskip4pt
\newpage

\nid{\bf 5. Conditional ZIPM = ZTP?}  

\nid Consider two subsets of the index set ${\cal K}$ and two subarrays of the data array ${\bf N}\equiv(N_{ij})$:
\begin{align*}
\Ome_Z^{\ne}&=\{(i,j)\,|\, Z_{ij}=1\},\\
\Ome_N^{\ne}&=\{(i,j)\,|\, N_{ij}\ne0\},\\
{\bf N}_Z^{\ne}&=(N_{ij}\,|\,Z_{ij}=1)=(M_{ij}\,|\,Z_{ij}=1),\\
{\bf N}^{\ne}&=(N_{ij}\,|\,N_{ij}\ne0)=(M_{ij}\,|\,M_{ij}\ne0).
\end{align*}
Both $\Ome_Z^{\ne}$ and $\Ome_N^{\ne}$ are random subsets, $\Ome_Z^{\ne}$ is unobserved, $\Ome_N^{\ne}$ is observed, and $\Ome_N^{\ne}\subseteq\Ome_Z^{\ne}$, so ${\bf N}^{\ne}\subseteq{\bf N}_Z^{\ne}$. Because ${\bf M}$ is independent of ${\bf Z}$, ${\bf N}_Z^{\ne}$ is a random subarray of the i.n.i.d. array $(M_{ij})$, where membership in this subarray depends only on ${\bf Z}$. Thus ${\bf N}^{\ne}$ is also is a (smaller) random subarray of the i.n.i.d. array $(M_{ij})$, where membership depends on both ${\bf Z}$ and the events $\{M_{ij}\ne0\}$.

The latter fact suggest a question: Is the conditional distribution of the two-component ZIPM rv
 $N_{ij}$ given $N_{ij}\ne0$ the same as the distribution  of the mixture of the conditional distributions of the two Poisson  components given that each is non-zero? The latter conditional distribution is the well-known {\it zero-truncated Poisson (ZTP) distribution}, also called positive Poisson, which has been thoroughly studied (cf. Johnson, Kemp, and Kotz (2005)). The ZTP distribution model also is an exponential family, with pmf given by
\begin{align}
g_\lam(x)=\frac{\lam^x}{(e^\lam-1)x!},\qquad x=1,2,\dots.\label{gpmf}
\end{align}

If the answer to the above question is yes, then estimation of $\pi,\mu,\nu$ and thus $\theta$ could be based on only the set of non-zero $N_{ij}$. That is, discard all 0's and view the remaining $N_{ij}$ as $\pi$-mixtures of two ZTP components with parameters $t_i\mu$ and $t_i\nu$. Because this involves only two mixture components rather than three as above, both being exponential families, and neither is degenerate, estimation methods such as the EM algorithm would be easier to carry out. 

Unfortunately the answer to the question is no. If we abbreviate $N_{ij}$ by $N$, $M_{ij}$ by $M$, and $Z_{ij}$ by $Z$, then the question can be exressed as follows: 
\begin{align*}
\mathrm{Is}&\ \ \P[N=x\,|\,N\ne0]=\frac{\pi\mu^x}{(e^\mu-1)x!}+\frac{(1-\pi)\nu^x}{(e^\nu-1)x!},\qquad x=1,2,\dots?
\end{align*}
However, for $x\ge1$,
\begin{align*}
P[N=x\,|\,N\ne0]&=\frac{\P[ZM=x,ZM\ne0]}{P[ZM\ne0]}\\
&=\frac{\P[M=x,Z=1,M\ne0]}{P[Z=1,M\ne0]}\\
&=\frac{\P[M=x,M\ne0]}{P[M\ne0]}\\
&=\frac{\frac{\pi\mu^x}{e^\mu x!}+\frac{(1-\pi)\nu^x}{e^\nu x!}} {1-\frac{\pi}{e^\mu}-\frac{(1-\pi)}{e^\nu}},
\end{align*}
since $M$ and $Z$ are independent, so the question becomes:
\begin{align*}
\mathrm{Is}&\ \ \frac{\frac{\pi\mu^x}{e^\mu x!}+\frac{(1-\pi)\nu^x}{e^\nu x!}} {1-\frac{\pi}{e^\mu}-\frac{(1-\pi)}{e^\nu}}=\frac{\pi\mu^x}{(e^\mu-1)x!}+\frac{(1-\pi)\nu^x}{(e^\nu-1)x!},\qquad x=1,2,\dots?
\end{align*}
After some algebra, this equation simplifies to
\begin{align*}
\Big(\frac{\mu}{\nu}\Big)^x&=\left(\frac{e^{\mu}-1} {e^\nu-1}\right),
\end{align*}
which cannot hold for all $x\ge1$ unless $\mu=\nu$.
\vskip10pt

 \nid{\it Acknowledgement.} I am grateful to Jon Wellner for his generous and always-insightful   comments.
 
\bigskip
\bigskip
\newpage

\centerline{\bf References}
\bigskip

\def\new{\hangindent=1.5pc}

\vskip .25cm \new
\noindent Arora, M.  and Chaganty, N. R. (2021). EM estimation for zero- and $k$-inflated Poisson regression model, {\it Computation}  {\bf 9} 94.

\vskip .25cm \new
\noindent Debavelaere, V.  and  S. Allassonni\`ere (2021). On the curved exponential family in the stochastic approximation expectation maximization algorithm, {\it ESAIM: Probability and Statistics} {\bf 25} 408-432.

\vskip .25cm \new
\noindent Efron, B.  and D. V. Hinkley (1978), Assessing the accuracy of the maximum likelihood estimator: observed versus expected Fisher information, {\it Biometrika}  {\bf 65} 457-482.

\vskip .25cm \new
\noindent Fahrmeir, L. (1987). Asymptotic likelihood inference for nonhomogeneous observations, {\it Statistische Hefte}  {\bf 28} 81-116.

\vskip .25cm \new
\noindent Guan, Y.  (2009). Variance stabilizing transformations of Poisson, binomial and negative binomial distributions, {\it Statist. Probability Letters} {\bf 79} 1621–1629.

\vskip .25cm \new
\noindent Hoadley, B. (1971). Asymptotic properties of maximum likelihood estimators for the independent not identically distributed case, {\it Ann. Math. Statist.}  {\bf 42} 1977-1991.

\vskip .25cm \new
\noindent Johnson, N. L.,  A. W. Kemp, and S. Kotz (2005). {\it Univariate Discrete Distributions (3rd ed.)}
], Wiley-Interscience, Hoboken, NJ.

\vskip .25cm \new
\noindent Lambert, D.  (1992). Zero-inflated Poisson regression, with an application to defects in manufacturing, {\it Technometrics} {\bf 34} 1-14.
 
\vskip .25cm \new
\noindent Laurent, S.  and C. Lagrand (2012). A Bayesian framework for the ratio of two Poisson rates in the context of vaccine efficacy trials, {\it ESAIM: Probability and Statistics} {\bf 16} 375-398.

\vskip .25cm \new
\noindent Li, H.-Q., M.-L. Tang, and W.-K. Wong (2014).  Confidence intervals for ratio of two Poisson rates using the method of variance estimates recovery, {\it Computational Statistics} {\bf 29} 869-889.

\newpage

\vskip .25cm \new
\noindent Louis, T. (1982). Finding the observed information matrix when using the EM Algorithm, {\it J. R. Statist.Soc. Series B} (1982) {\bf 44} 226-233.

\vskip .25cm \new
\noindent McLachlan, G. J. and T. Krishnan (2008). {\it The EM Algorithms and its Extensions, 2nd ed.}, New York: Wiley.

\vskip .25cm \new
\noindent Robert, C. and G. Casella (2004). {\it Monte Carlo Statistical Methods}, New York: Springer-Verlag.

\vskip .25cm \new
\noindent Silverman, B. W. (1986). {\it Density Estimation for Statistics and Data Analysis}, London: Chapman \& Hall/CRC.

\end{document}